\documentclass[12pt,reqno]{article}

\usepackage[usenames]{color}
\usepackage{amssymb}
\usepackage{amsmath}
\usepackage{amsthm}
\usepackage{amsfonts}
\usepackage{amscd}
\usepackage{graphicx}

\usepackage[colorlinks=true,
linkcolor=webgreen,
filecolor=webbrown,
citecolor=webgreen]{hyperref}

\definecolor{webgreen}{rgb}{0,.5,0}
\definecolor{webbrown}{rgb}{.6,0,0}

\usepackage{color}
\usepackage{fullpage}
\usepackage{float}

\usepackage{graphics}
\usepackage{latexsym}
\usepackage{epsf}
\usepackage{breakurl}

\newcommand{\seqnum}[1]{\href{https://oeis.org/#1}{\rm \underline{#1}}}
\newcommand{\C}{{\mathcal C}}
\newcommand{\N}{{\mathbb N}}
\newcommand{\Prime}{{\mathcal P}}

\begin{document}


\theoremstyle{plain}
\newtheorem{theorem}{Theorem}
\newtheorem{corollary}[theorem]{Corollary}
\newtheorem{lemma}[theorem]{Lemma}
\newtheorem{proposition}[theorem]{Proposition}

\theoremstyle{definition}
\newtheorem{definition}[theorem]{Definition}
\newtheorem{example}[theorem]{Example}
\newtheorem{conjecture}[theorem]{Conjecture}

\theoremstyle{remark}
\newtheorem{remark}[theorem]{Remark}

\begin{center}
\vskip 1cm{\LARGE\bf Conjectures about Primes and Cyclic Numbers
}
\vskip 1cm
\large
Joel E. Cohen\footnote{Author's other affiliations: Departments of Statistics, Columbia University \& University of Chicago.} \\
The Rockefeller University \\
1230 York Avenue, Box 20\\
New York, NY 10065\\
USA \\
\href{mailto:email}{\tt cohen@rockefeller.edu} \\
\end{center}

\vskip .2 in
\begin{abstract}
 A positive integer $n$ is defined to be \emph{cyclic} if and only if 
every group of size $n$ is cyclic.
Equivalently, $n$ is cyclic if and only if
$n$ is relatively prime to the number of positive 
integers less than $n$ that are relatively prime to $n$.  
Because every prime number is cyclic,
it is natural to ask whether a (proved or conjectured) property of primes 
extends to cyclic numbers.
I review  proved or conjectured properties of primes 
(including some new conjectures about primes) and
propose analogous conjectures about cyclic numbers.
Using the 28,488,167 cyclic numbers less than $10^8$,
I test the conjectures about cyclic numbers and disprove 
the cyclic analog of the second conjecture about primes of Hardy and Littlewood.
Proofs or disproofs of the remaining conjectures are invited.
\end{abstract}

\section{Introduction}
I propose some conjectures about cyclic numbers $\C:=(c_1,\ c_2,\ \ldots )$ (sequence 
\cite[\seqnum{A003277}]{oeis} in the \emph{On-Line Encyclopedia of Integer Sequences} (OEIS))
based on analogous proved or conjectured properties of prime numbers $\Prime:=(p_1,\ p_2,\ \ldots )$
(\seqnum{A000040}).
I test the conjectures about cyclic numbers (or, for brevity, cyclics)
using the 28,488,167 cyclics less than $10^8$.
I also test some new conjectures about prime numbers (or, for brevity, primes) using the 50,847,534 primes less than $10^9$.
I invite proofs, disproofs, further numerical confirmations, counterexamples,
and additional conjectures.

A natural number (positive integer)
$n\in\N:=\{1, 2, 3, \ldots\}$ is \emph{cyclic} if and only if
there exists only one group of size $n$, up to isomorphism. 
If  $\gcd$ is the greatest common divisor and 
$\varphi(n)$ is Euler's totient function (the number of positive 
 integers less than $n$ that are relatively prime to $n$, \seqnum{A000010}), then, according to Tibor Szele (1918--1955) \cite{Szele1947},
the number
$n\in\N$ is cyclic if and only if 
$\gcd(n,\varphi(n)) = 1$.

I use Szele's condition to compute which $n\in\N$ are cyclic numbers.
My first $10^4$ computed cyclic numbers exactly match
the first $10^4$ cyclic numbers computed independently by T. D. Noe 
in \seqnum{A003277}.

Michel Lagneau \cite[\seqnum{A003277}, November 18, 2012]{oeis} asserted without proof that
$n\in\N$ is {cyclic} if and only if
$\varphi(n)^{\varphi(n)} \equiv 1$ mod $n$. 
Richard P. Stanley (personal communication, 2025-01-20 20:00) gave an elegant short proof of Lagneau's condition. 
I quote it with his permission.
First, assume that $n$ is not cyclic, so $\gcd(n,\varphi(n))=d>1$. 
Then $\varphi(n)^{\varphi(n)}$ is divisible by $d$ so cannot be 
congruent to 1 mod $n$. 
On the other hand, assume that $n$ is cyclic. 
Euler's generalization of Fermat's little theorem 
implies that $k^{\varphi(n)}\equiv 1$ mod $n$ whenever $\gcd(k,n)=1$. 
Putting $k=\varphi(n)$ completes the proof.

Determining which numbers are cyclic 
by evaluating $\varphi(n)^{\varphi(n)}$ by brute force can be problematic. 
If $n$ is a large prime, then 
$\varphi(n)^{\varphi(n)} = (n-1)^{n-1}$ is very large.
Such a large number is likely to be imprecisely approximated by 
software that does not offer infinite precision arithmetic for integers, 
and is likely to be very cumbersome or impossible for software that does offer infinite precision arithmetic for integers. 
For example, $p_{5761455} = 99,999,989 = c_{28488165}$, 
and $\varphi(99999989)=99,999,988$ differs from $10^8$ by roughly one part in a million.
To a first approximation, $99999988^{99999988}$ is an integer with 
roughly $1 + 8\times 10^8$ decimal digits.
Matlab \cite{matlab} approximates $99999988^{99999988}$ as infinity.

Alexei Kourbatov (personal communications, 2025-05-24 09:18, 2025-06-01 13:46)
pointed out that 
Lagneau's condition can be evaluated using an algorithm \cite[p.\ 71, algorithm 2.143]{HAC1996} for modular 
exponentiation \cite{WikiModularExponentiation} that
requires no intermediate results of order $O(n^{n-1})$.
Whether using Szele's criterion (as I do here) or Lagneau's criterion
for a number to be cyclic, the most computationally expensive step 
is finding $\varphi(n)$.


The sequence $\C$ of cyclics begins 
$(1, 2, 3, 5, 7, 11, 13,15, 17, 19, 23, 29, 31, 33, 35, 37,\ldots)$.
Cyclics $\C$ are the union of primes $\Prime$
and the composite numbers $n\in\N$ such that $n$ and $\varphi(n)$
are relatively prime or coprime (\seqnum{A050384}, e.g., 
1, 15, 33, 35, 51, 65, 69, 77, 85, 87, 91, 95, 115, 119, 123, 
 133, 141, 143, 145, 159,  $\ldots$).
The only cyclic that is a square is $c_1=1$.
The only cyclic that is even is $c_2=2$.
Consequently, the only cyclic of the form $n(n-1)$ for $n\in\N$ is $c_2 = 2$ with $n=2$, and
no cyclic is of the form $n(n+1)$ for $n\in\N$, because both $n(n-1)$ and $n(n+1)$ are even.

For an increasing integer sequence $a:=(a(1),\ a(2),\ a(3),\ \ldots)$,
the counting function of $a$ evaluated at a positive real number $x$ is the number
of elements of $a$ that are less than or equal to $x$.
For example, the counting function $\pi(\cdot)$ of primes $\Prime$ satisfies $\pi(10)=4$.
The prime number theorem \cite{Hadamard1896, delaValleePoussin189900} gives 
that $\pi(x) \sim x/\log x$ as $x\to\infty$.
Let 
\begin{equation}\label{eq:cycliccountingfn}
C(x):= \sum_{\substack{ m \leq x \\ \text{$m$ cyclic} }} 1
\end{equation}
be the counting function of cyclic numbers, that is, the number of
cyclic numbers that do not exceed positive real $x$ (\seqnum{A061091}).
For $n=1,\ldots,20$, $C(n)= 1, 2,  3,     3,     4,     4,     5,     5,     5,     
5,     6,     6$,     $7,     7,     8,     8,     9,     9,    10,    10$. 
P\'al Erd\H{o}s (1913--1996) \cite{Erdos1948} proved that 
\begin{equation}\label{eq:ErdosCountingFn}
C(x) \sim \frac{x}{e^{\gamma}\log \log \log x} \textrm{ as } x\to\infty.
\end{equation}
Here $\gamma\approx 0.5772156649\ldots$ is the Euler--Mascheroni constant and
$e^{\gamma}\approx 1.78107241799\ldots$ . 
Paul Pollack \cite{Pollack2022} gave an asymptotic series expansion 
\begin{equation}\label{eq:PollackCountingFn}
C(x)  \sim \frac{x}{e^{\gamma}\log \log \log x}\left(1 - \frac{{\gamma}}{\log \log \log x} +
\cdots
\right) \textrm{ as } x\to\infty
\end{equation}
with additional terms. I shall use just these first two terms.

Cyclics are much more abundant than primes asymptotically because
$\lim_{x\to\infty} \pi(x)/C(x) = 0$.
Hence, asymptotically, almost all cyclics are composite.
John Campbell and I \cite{CampbellCohen2025} observed that,
 since $C(c_n) = n$ by definition,
\eqref{eq:ErdosCountingFn} implies that 
\begin{align}
 &  c_n  \sim e^{\gamma}n\log \log \log n \textrm{ as } n\to\infty, \label{eq:cn} \\ 
 & \lim_{n\to\infty} \frac{c_{n+1}}{c_n} = 1, \nonumber \\ 
 & \lim_{x\to\infty}\frac{x}{c_{C(x)} } =1,\label{eq:xoncC}\\ 
 & \lim_{x\to\infty}\frac{\log C(x)}{\log x} =1. \nonumber 
\end{align}

The observation that cyclics are asymptotically much more abundant than primes 
motivates investigating which (proved or conjectured) properties of primes depend
on their asymptotic scarcity relative to cyclics, 
and which properties of primes carry over 
(exactly or asymptotically) to the  more abundant cyclics.

Another infinite increasing integer sequence that contains all primes is $\N$.
But $\N$ does not share an elementary property that $\Prime$ and $\C$ share,
namely, that the only even element of the sequence is 2.
Similarly, while $\Prime$ contains no squares and $\C$ contains exactly one 
square, $\N$ includes infinitely many squares.
Other infinite increasing integer sequences that share important properties with 
$\Prime$ and $\C$ remain to be investigated.

Campbell and I \cite{CampbellCohen2025} proved two 
analogies between primes and cyclics.
First, under the Bernhard Riemann (1826--1866) hypothesis, 
the $n$th prime gap satisfies 
$p_{n+1} - p_{n} = O(\sqrt{p_{n}} \log p_{n})$ as $n\to\infty$ \cite{Cramer1921}.
More precisely, under the Riemann hypothesis, for every $p_n>3$,
$p_{n+1} - p_{n} < \frac{22}{25}\sqrt{p_n}\log p_n$ \cite{Carneiro2019}.
We \cite[Theorem 2]{CampbellCohen2025} proved that, under the Riemann hypothesis, 
the first difference of consecutive cyclics satisfies
$c_{n+1} - c_{n} = o\left( \sqrt{p_n}  \log p_n \right)$.
Second, if $m_n(\Prime)$ is the mean and $v_n(\Prime)$ is the variance of the first $n$ primes, then asymptotically
$v_n(\Prime) \sim (1/3)(m_n(\Prime))^2$ as $n\to\infty$
\cite{Cohen2016}.
We \cite[Theorem 1]{CampbellCohen2025} proved, 
without the Riemann hypothesis, that the mean $m_n(\C)$
and the variance $v_n(\C)$
of the first $n$ cyclics satisfy the same asymptotic relationship,
$v_n(\C) \sim (1/3)(m_n(\C))^2$ as $n\to\infty$.

This project of generalizing from primes to cyclics is not guaranteed to succeed.
After proposing in section \ref{sec:conjectures4cyclics} conjectures that numerical calculations have so far failed to reject,
I give in section \ref{sec:HLconjecture2} a counterexample to show that the
analog for cyclics of the second conjecture \cite{HardyLittlewood1923} of 
Godfrey Harold Hardy (1877--1947) and John Edensor Littlewood (1885--1977)
fails.
This counterexample provides a further small, indirect hint in support of the belief of Hensley and Richards
\cite{HensleyRichards1974}
that the second conjecture of Hardy and Littlewood for primes is  false.

Because $\Prime \subset \C$,
if infinitely many primes have property X, then 
infinitely many cyclics have property X.
For example, Euler's proof that $\lim_{n\to\infty}\sum_{j=1}^n p_j^{-1} = \infty$
immediately implies that $\lim_{n\to\infty}\sum_{j=1}^n c_j^{-1} = \infty$.
But if every prime has property X, it 
may be true or false, depending on property X,
that every cyclic has property X.
Conversely, if infinitely many cyclics have property X, it 
may be true or false, depending on property X,
that infinitely many primes have property X.
But if every cyclic has property X, then
every prime has property X.

Consequently, when a conjecture about primes has been extensively verified
numerically, if
that conjecture immediately implies
the corresponding conjecture about cyclics,
there is no need, and I do not bother, to test numerically 
the analogous conjecture about cyclics. 
I test numerically only those conjectures about cyclics
not immediately implied by properties of primes
that are proved or conjectured and numerically supported.

A helpful referee pointed out that many additional questions could be asked about cyclics.
For example, the referee asked,
are the cyclics equidistributed over arithmetic progressions of a prescribed modulus?
How does the sum of all cyclics less than or equal to positive real $x$ behave as a function of $x$?
The latter question leads to the first and only theorem of this paper, 
which reports, for a fixed positive integer $k$, 
the sum of the $k$th power of all cyclics less than or equal to a positive real $x$ as a function of $x$.

\begin{theorem}\label{thm:sumofcyclics}
Fix $k\in\N$. For $n\in\N$, if $c_1,\ldots,c_n\in\C$ are the first $n$ cyclic numbers, then
\begin{equation}\label{eq:sumofcyclics}
c^k_1+\cdots+c^k_n \sim \frac{nc^k_n}{k+1} \sim 
 \frac{n^{k+1}e^{k\gamma}(\log \log \log n)^k}{k+1} \textrm{ as } x\to\infty.
\end{equation}
In particular, $c_1+\cdots+c_n \sim nc_n/2\sim n^2e^{\gamma}\log \log \log n/2$.

For positive real $x$, as $x\to\infty$, 
the sum of the $k$th power of all cyclics less than or equal to $x$ is
asymptotic to
\begin{equation}\label{eq:sumofcyclicstopowerkLTx}
\frac{C(x)c^k_{C(x)}}{k+1}\sim \frac{C(x)x^k}{k+1}
 \sim \frac{x^{k+1}}{(k+1)e^{\gamma}\log \log \log x}\left(1 - \frac{{\gamma}}{\log \log \log x}
\right),
\end{equation}
using \eqref{eq:PollackCountingFn}.
In particular, the sum of all cyclics less than or equal to positive real $x$ 
is asymptotic to
\begin{equation}\label{eq:sumofcyclicsLTx}
 \frac{C(x)x}{2}
 \sim \frac{x^{2}}{2e^{\gamma}\log \log \log x}\left(1 - \frac{{\gamma}}{\log \log \log x}
\right).
\end{equation}
\end{theorem}
\begin{proof}
Campbell and I \cite[Theorem 1]{CampbellCohen2025} 
showed that 
\begin{equation*}
n^{-1}(c_1^k+\cdots+c_n^k) \sim\frac{c_n^k }{k+1} \textrm{ as } x\to\infty.
\end{equation*}
Hence,  using \eqref{eq:cn},
$c^k_1+\cdots+c^k_n \sim nc_n^k /(k+1)
\sim n( e^{\gamma}n\log \log \log n)^k /(k+1)$\\
$=
n^{k+1}e^{k\gamma}(\log \log \log n)^k/(k+1)$, 
proving \eqref{eq:sumofcyclics}.

Replacing $n$ in \eqref{eq:sumofcyclics} by $C(x)$ and 
using \eqref{eq:xoncC} to approximate $c_{C(x)}$ yields 
\eqref{eq:sumofcyclicsLTx}.
\end{proof}

\begin{figure}[htb]
\centering
\includegraphics[width=1\textwidth]{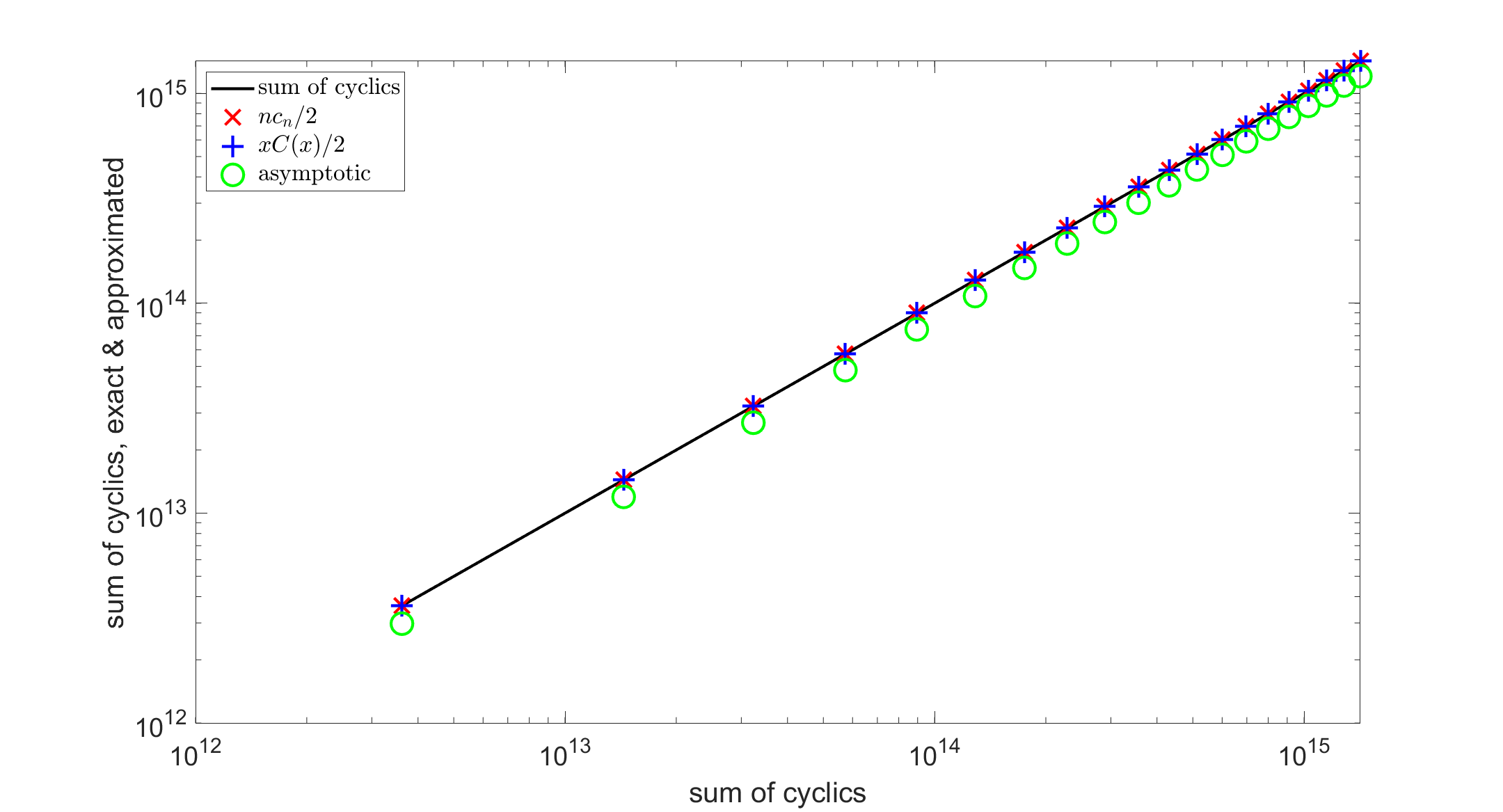}
\caption{For $m = 1, 2,\ldots, 20$ and $x=m\times 5 \times 10^6$,
the abscissa of each plotted point is the exact sum of cyclics less than or equal to $x$.
The ordinate of each point compares the exact sum of cyclics less than or equal to $x$
(solid black line, with abscissa = ordinate) with three approximations:
$nc_n/2$ (red $\times$ marker); $xC(x)/2$ (blue + marker); and 
the asymptotic approximation on the right side of \eqref{eq:sumofcyclicsLTx} (open green circle).}
\label{fig:SumOfCyclicsLessThan_x}
\end{figure}

In general, $xC(x)/2 > nc_n/2$ because in general $x>c_n$ while $C(x)=n$.
Figure \ref{fig:SumOfCyclicsLessThan_x} shows
that $xC(x)/2$ and $nc_n/2$ closely approximate the exact sum of cyclics, and
the asymptotic approximation on the right side of \eqref{eq:sumofcyclicsLTx}
consistently falls below
the exact sum of cyclics and the other two approximations.

\section{Conjectures}\label{sec:conjectures4cyclics}
\subsection{Landau's list and Legendre's relatives}

In 1912, Edmund Landau (1877--1938) \cite{Landau1912} presented four historic conjectures about primes.
Extensively verified numerically, these conjectures are generally believed to be true but unproved as of 2025 \cite{WikiLandau}. 
Many unconfirmed and unrefuted claims to have proved one or several of Landau's conjectures
have been published or posted but I am not aware that such a claimed
proof has been independently confirmed.

If true, each of Landau's four conjectures about primes  immediately implies the  conjecture 
about cyclics that follows it below. 
But the conjecture about cyclics may be true even if the corresponding 
Landau conjecture about primes is false.
Here are Landau's four conjectures about primes 
and analogous conjectures about cyclics.
I include a few novel conjectures about cyclics suggested by these
analogs of Landau's conjectures.

\subsubsection{Landau's problem 1}

First, Christian Goldbach (1690--1764) conjectured that every even $n\in\N$ greater than 2 is a 
sum of two primes.

\begin{conjecture}[Goldbach analog for cyclics]\label{conj:Goldbach}
Every even $n\in\N$ greater than 2 is a sum of two cyclics.
\end{conjecture}

On seeing Conjecture \ref{conj:Goldbach} in a prior draft of this paper, 
Carl Pomerance (personal communication, 2025-03-05) \cite{Pomerance2025}
proved that every sufficiently large even $n$ is a sum of two cyclics. 
I quote his result with his permission.
Let $G(n)$ be the number of pairs of cyclics $c_1,\ c_2$ such that $c_1+c_2 = n$.
Then if $n$ is even, Pomerance \cite{Pomerance2025} proved,
$$G(n) \sim \prod_{p>2}\left(1-\frac{1}{(p-1)^{2}}\right)\cdot
\frac{2n}{(e^{\gamma}\log \log \log n)^{2}}\cdot
\prod_{\substack{p\mid n \\2<p<\log \log n}}\frac{p-1}{p-2}
 \text{ as }n\to\infty.$$

\subsubsection{Landau's problem 2}

Second in Landau's list, the twin prime conjecture states that there are infinitely many primes $p$ 
such that $p+2$ is also a prime.
Here $p$ and $p+2$ are called twin primes.

\begin{conjecture}[twin cyclics analog]\label{conj:twincyclics}
There are infinitely many cyclics $c\in\C$ such that $c+2$ is also cyclic.
\end{conjecture}

On seeing Conjecture \ref{conj:twincyclics} in a prior draft of this paper, 
Carl Pomerance (personal communication, 2025-03-05) \cite{Pomerance2025}
proved a much stronger result, which I quote with his permission.
For positive real $x$, let $C_2(x)$ be the number of cyclics $c\leq x$
such that $c+2$ is also cyclic. Then as $x\to\infty$,
\begin{equation}\label{PomeranceA2}
C_2(x) \sim 2\prod_{p>2}\left(1-(p-1)^{-2}\right)x(e^{\gamma}\log \log \log x)^{-2}.
\end{equation}
The right side of \eqref{PomeranceA2} approaches infinity as $x\to\infty$, proving Conjecture \ref{conj:twincyclics}.

It is well known that the only prime $p$ such that $p,\ p+2,\ p+4$ are all primes is $p=3$,
because if $p>3$, one of $p,\ p+2,\ p+4$ must be divisible by 3.
Cyclics are different. The composite cyclics (\seqnum{A050384}) include multiple triplets, such as
141, 143, 145, and 213, 215, 217, and 319, 321, 323, and 391, 393, 395.

\begin{conjecture}[cyclic triplets]\label{conj:cyclictriplets}
There are infinitely many composite cyclics $c\in\C$ such that $c,\ c+2,\ c+4$ are all composite cyclic.
\end{conjecture}
Carl Pomerance (personal communication, 2025-03-05) \cite{Pomerance2025}
proved a related result, which I quote with his permission: there are infinitely many cyclic triplets $c,\ c+2,\ c+4$ 
(not necessarily all composite cyclics, as in Conjecture \ref{conj:cyclictriplets}),
and their counting function is of order $x(\log \log \log x)^{-5/2}(\log x)^{-1/2}$ as $x\to\infty$. 
As Carl Pomerance pointed out (personal communication, 2025-05-24 10:43), 
since the number of primes up to $x$ is of order $x/\log x$,
which is much smaller, asymptotically most of these triplets consist of three
composites.

Further, the cyclics (whether prime or composite) include multiple quintuplets with successive gaps equal to 2, such as 11, 13, 15, 17, 19, and 29, 31, 33, 35, 37, and
65, 67, 69, 71, 73, and 83, 85, 87, 89, 91, and 137, 139, 141, 143, 145, and 209, 211, 213, 215, 217, and 263, 265, 267, 269, 271.

\begin{conjecture}[cyclic quintuplets]\label{conj:cyclicquintets}
There are infinitely many cyclics $c\in\C$ such that $c,\ c+2,\ c+4,\ c+6,\ c+8$ are all cyclic.
\end{conjecture}

Carl Pomerance (personal communication, 2025-03-05) \cite{Pomerance2025}
proved a general result which implies, 
for example, that infinitely many all-cyclic 8-tuples have the form
$n, n+2, n+4, n+6, n+8, n+10, n+12, n+14$ but there are no all-cyclic 9-tuples
with the additional term $n + 16$ because one of these nine
numbers is divisible by 9, therefore not cyclic.

\subsubsection{Landau's problem 3}

Third in Landau's list, Adrien-Marie Legendre (1752--1833) \cite{Legendre1808} conjectured that for every $n\in \N$, there exists a prime $p\in\Prime$
such that $n^2<p<(n+1)^2$. 
Legendre's claimed proof was based on a prior claim that was false.
Several claims to prove Legendre's conjecture have been published or posted
but I am not aware that any has been independently confirmed.

\begin{conjecture}[Legendre analog for cyclics]\label{conj:Legendre}
For every $n\in \N$, there exists a cyclic $c\in\C$
such that $n^2<c<(n+1)^2$.
\end{conjecture}
A cyclic $c$ such that $n^2<c<(n+1)^2$ exists for all $1 \leq n\leq 9998$. 
The use of strict inequality in Conjecture \ref{conj:Legendre} is justified because
no cyclic other than 1 is a square.

Honor\'e Adolphe Desboves \cite[Theorem 2, p.\ 290]{Desboves1855} 
(1818--1888) conjectured that
for every $n\in \N$, there exist two primes $p, p'$ 
such that $n^2<p < p'<(n+1)^2$.
Desboves asserted that if Legendre's conjecture is true, then his conjecture follows.
The converse is obvious.
I confirmed Desboves' conjecture numerically for the 50,847,534 primes less than $10^9$.
If true, Desboves' conjecture immediately implies Conjecture \ref{conj:Desbovescyclic} about cyclics. 
But Conjecture \ref{conj:Desbovescyclic} may be true even if Desboves' conjecture is false.

\begin{conjecture}[Desboves analog for cyclics]\label{conj:Desbovescyclic}
For every $n\in \N$, there exist two cyclics $c, c'\in\C$ 
such that $n^2<c < c'<(n+1)^2$.
\end{conjecture}
For every positive integer $n\leq 3161$, there exist $c, c'\in \C$ 
such that $n^2<c < c'<(n+1)^2$.

The next  conjectures generalize Legendre's and Desboves's conjectures to $k>2$ primes and cyclics.

For $n\in\N$, let $N_{\Prime}(n):= \#\{p\in\Prime \mid p\in(n^2,(n+1)^2)\}$ be
the number of primes in the interval $(n^2,(n+1)^2)$ (\seqnum{A014085} apart from an initial 0).
For example, in the first 25 intervals $(n^2,(n+1)^2),\ n=1,\ldots,25$,
the numbers of primes are, respectively, 
2,    2,    2,    3,    2,    4,    3,    4,    3,    5,    4,    5,    5,    4,    6,
7,    5,    6,    6,    7,    7,    7,    6,    9,    8.
For example, $N_{\Prime}(6) = 4$ because four primes,
37, 41, 43, 47, are between $6^2 = 36$ and $7^2 = 49$.

In a prior draft of this paper, I conjectured that
$N_{\Prime}(n)$ is asymptotic (as $n\to\infty$) to a regularly varying function 
of $n$ with positive index not exceeding 1.
Recall that a regularly varying function \cite{Feller1971, Seneta1976, Kevei2019}
maps the positive half line $x>0$ into
the positive half line and takes the form  $x\mapsto x^{\rho}\ell(x)$.
The exponent $\rho$ of $x$ is a real number, 
commonly called the index of the regularly varying function, and 
$\ell(x)$ is a slowly varying function of $x$, that is, for every $\lambda>0$,
$\lim_{x\to\infty}\ell(\lambda x)/\ell(x) = 1$.
A regularly varying function generalizes a power function $x\mapsto x^{\rho}$.

On seeing this conjecture, 
Pierre Deligne (personal communication, 2025-03-06 19:12) refined my conjecture to
a much more specific, heuristically plausible conjecture,
which I quote with his permission.
Deligne observed that the length of the interval $(n^2,(n+1)^2)$
is asymptotic to $2n$
and the probability that an integer in this interval is prime is
asymptotic to $1/\log(n^2)$, so the number of primes in 
$(n^2,(n+1)^2)$ should be asymptotic to 
$2n/\log(n^2) = n/\log n$.
Deligne commented, ``Of course with no way to prove it.''

\begin{conjecture}[Deligne's conjecture: primes in intervals between successive squares]\label{conj:Np}
As $n\to\infty$,
the number $N_{\Prime}(n)$ of primes in the interval $(n^2,(n+1)^2)$
satisfies $N_{\Prime}(n)\sim n/\log n$.
\end{conjecture}

\begin{figure}[htb]
\centering
\includegraphics[width=1\textwidth]{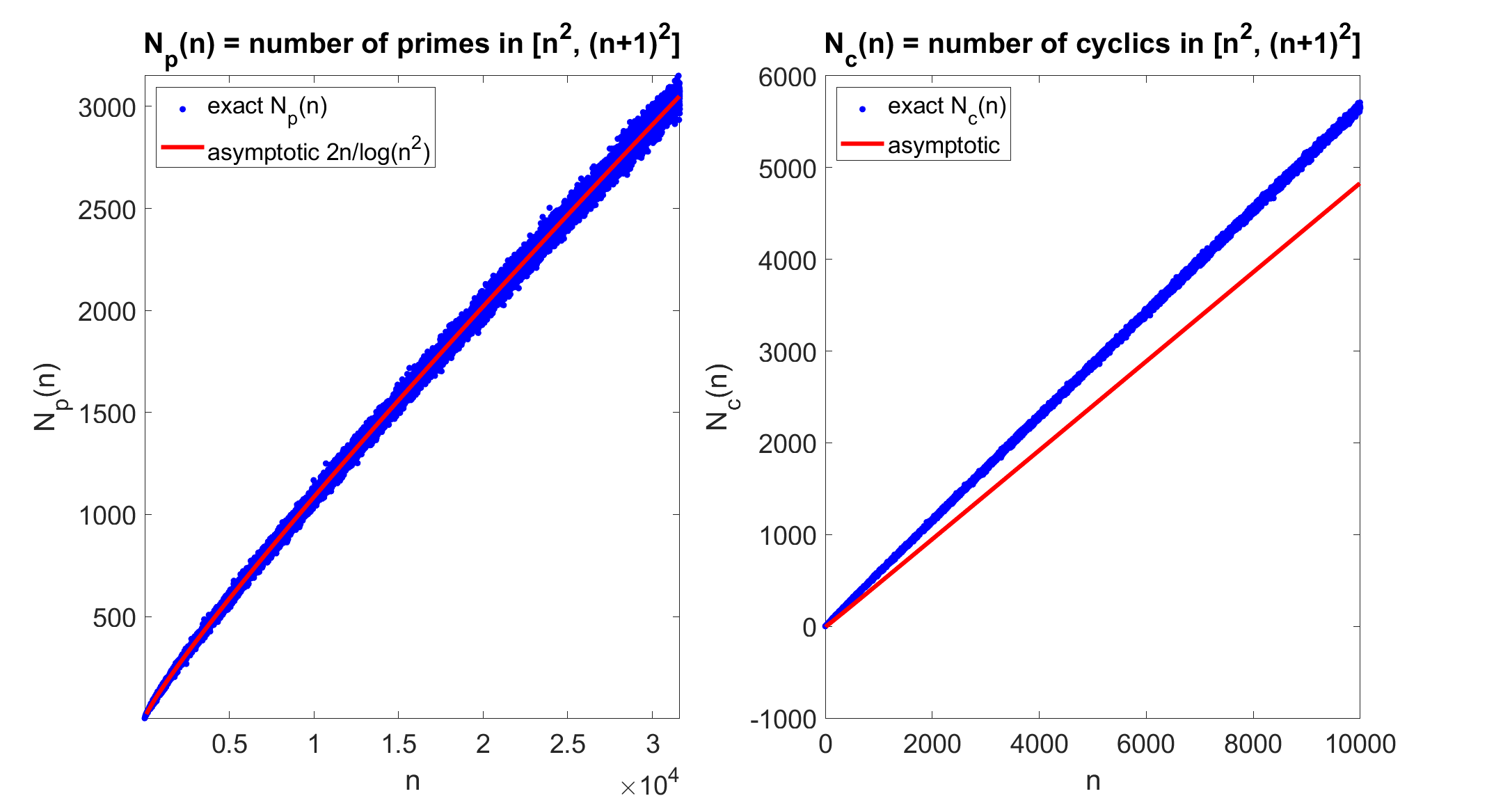}
\caption{(left) The number $N_{\Prime}(n)$ of primes in the interval $(n^2, (n+1)^2)$ 
(blue dots) for $n= 1,\ldots, 31 621$ and a conjectured asymptotic approximation $n/\log(n)$ (red curve).
(right) The number $N_c(n)$ of cyclics in the interval $(n^2, (n+1)^2)$
(blue dots) for $n= 1,\ldots, 9998$, and a conjectured asymptotic approximation \eqref{eq:CyclicsInSquaresUsingPollack} (red line).}
\label{fig:PrimeCyclicPowerLaws}
\end{figure}

Figure \ref{fig:PrimeCyclicPowerLaws} (left)
plots $N_{\Prime}(n)$ for $n= 1,\ldots, 31 621$ (blue dots)
and the asymptotic approximation $n/\log(n)$ (red line).
The results support Deligne's Conjecture \ref{conj:Np} for primes.

Imitating Deligne's heuristic argument for primes, 
the length of the interval $(n^2,(n+1)^2)$
is asymptotic to $2n$.
The probability that an integer in this interval is cyclic should
be proportional to $C((n+1)^2)/(n+1)^2$, which by 
\eqref{eq:PollackCountingFn} is asymptotic to
\begin{equation*}
 \frac{1}{e^{\gamma}\log \log \log ((n+1)^2)}\left(1 - \frac{{\gamma}}{\log \log \log ((n+1)^2)} \right).
\end{equation*}
Hence as $n\to\infty$, the number of cyclics in 
$(n^2,(n+1)^2)$ should be asymptotic  to 
\begin{equation}\label{eq:CyclicsInSquaresUsingPollack}
 \frac{2n}{e^{\gamma}\log \log \log ((n+1)^2)}\left(1 - \frac{{\gamma}}{\log \log \log ((n+1)^2)} \right)
\sim  \frac{2n}{e^{\gamma}\log \log \log (n)}\left(1 - \frac{{\gamma}}{\log \log \log (n)} \right).
\end{equation}
This heuristic argument is not a proof of the following conjecture.

\begin{conjecture}[cyclics in intervals between successive squares]\label{conj:Nc}
As $n\to\infty$, the number $N_c(n):= \#\{c\in\C\mid c\in(n^2,(n+1)^2)\}$ of cyclics in the interval $(n^2,(n+1)^2)$
is asymptotic to \eqref{eq:CyclicsInSquaresUsingPollack}.
\end{conjecture}

In my numerical calculations, the ratio of $N_c(n)$ to the 
corresponding function of $n$ in \eqref{eq:CyclicsInSquaresUsingPollack}
declines slowly toward 1 as $n$ increases, 
despite the increase in Figure \ref{fig:PrimeCyclicPowerLaws} (right) 
in the arithmetic difference between $N_c(n)$ 
and the corresponding function of $n$ in \eqref{eq:CyclicsInSquaresUsingPollack}.

Let $L_p$ be the OEIS sequence \seqnum{A349997}, defined as
``Numbers $k$ such that the number of primes in any [i.e., every] interval $[j^2,(j+1)^2], j>k$, is not less than the number of primes in the interval $[k^2,(k+1)^2]$.'' 
Let $L_p(n)$ be the $n$th element of $L_p$ in increasing order.

\begin{conjecture}[$k$-fold Legendre for primes]\label{conj:Legendrek1}
For primes,
$L_p =\{1, 7, 11, 17, 18, 26, 27, 32$, $46, 50$, $56, 58, 85, 88, 92, 137, 143, 145,\ldots \}.$
\end{conjecture}

For example, $L_p(2) = 7$ means that $N_p(7) = 3$
(i.e., three primes 53, 59, 61 lie between $7^2$ and $8^2$)
and (conjecturally, based on available computations) for every $j>7$, $N_p(j)\geq 3$.

\begin{figure}[htb]
\centering
\includegraphics[width=1\textwidth]{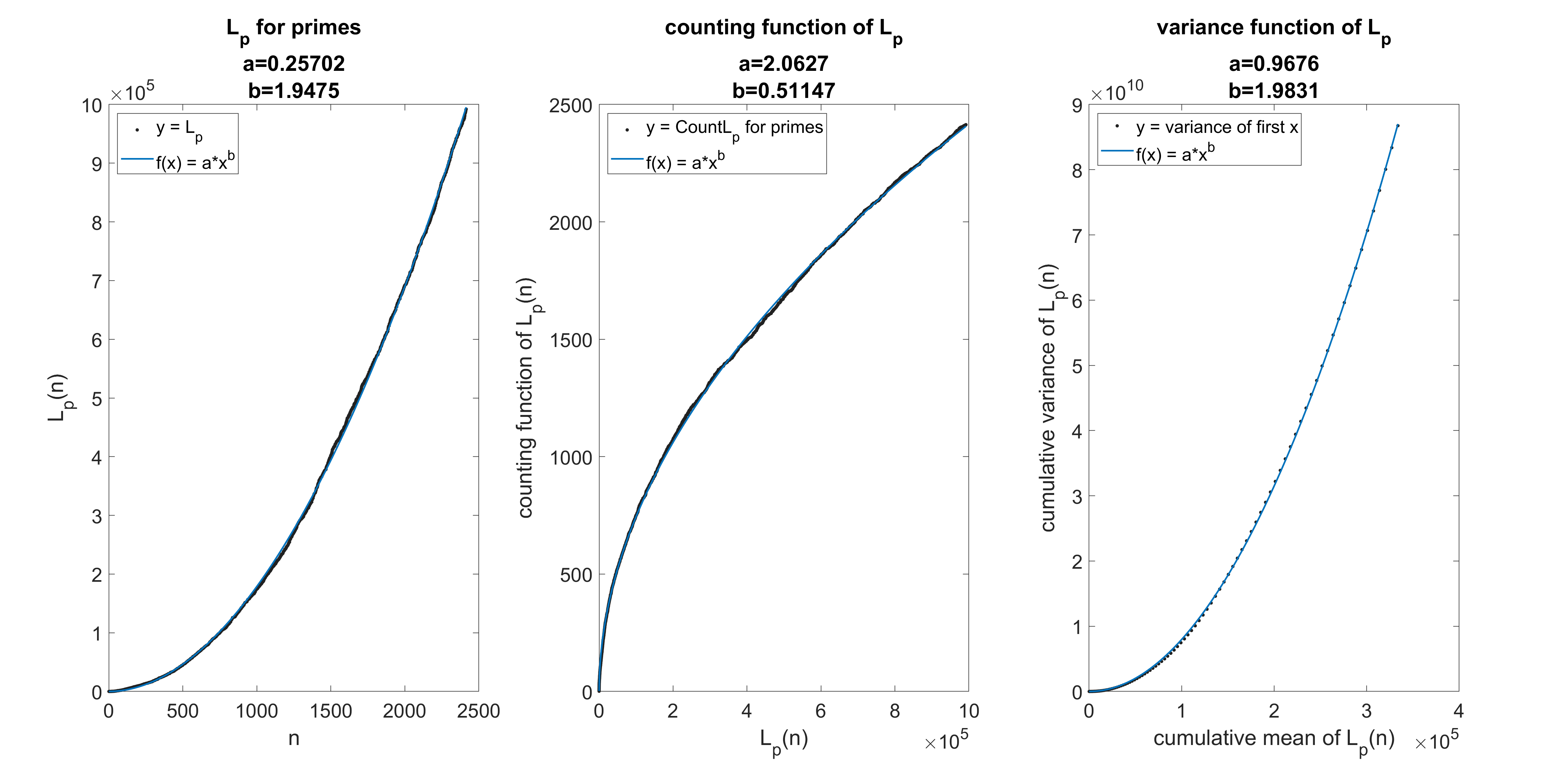}
\caption{(left) The sequence $L_p$ (\seqnum{A349997}, black dots), and 
the power function (blue curve) fitted by least squares.
(middle) The counting function of $L_p$ (black dots) and
the power function (blue curve) fitted by least squares.
(right) Variance function of $L_p$ (black dots) and a power function (blue curve) 
fitted by least squares.}
\label{Lprimes}
\end{figure}

Hugo Pfoertner tabulated 2414 (conjectural) values of $L_p$ at OEIS \seqnum{A349997.b349997.txt}
without reporting the number of primes he considered. 
These numerical values are conjectural because they depend
on an infinite sequence of primes not 
accessible to computation and not yet analyzed mathematically.
Pfoertner's 2414 values 
appear (Figure \ref{Lprimes} left) to be well approximated by $an^b$ with 
$a\approx 0.257,\ b \approx 1.9475$.

\begin{conjecture}[asymptotic $k$-fold Legendre for primes]\label{conj:Legendrekasymptotic}
As $n\to\infty$, $L_p(n)$ 
is asymptotic to a regularly varying function 
with index $b$ that satisfies $3/2 < b \leq 2$.
\end{conjecture}

The counting function of the sequence $L_p(\cdot)$ is 
defined for each $m\in\N$ as $\sum_{L_p(n)\leq m}1.$  
Because $L_p(\cdot)$ and its counting function are asymptotically inverses,
a mathematical consequence of Conjecture \ref{conj:Legendrekasymptotic} is that
the counting function of $L_p(\cdot)$ is asymptotic to a regularly varying function
with index $1/2 \leq 1/b < 2/3$ 
 \cite[pp.\ 21--27]{Seneta1976} \cite[section 8, pp.\ 16--17]{Kevei2019}.
Define $m_p(n)$ and $v_p(n)$ to be, respectively, 
the mean and the variance of $L_p(1), L_p(2),\ldots,L_p(n)$.
Then applying \cite[Theorem 1]{Cohen2023}
to Conjecture \ref{conj:Legendrekasymptotic} gives
$v_p(n) \sim m_p(n)^2/\bigl((1/b)(1/b+2)\bigr) = b^2m_p(n)^2/(1+2b)$.

The estimated exponent $b\approx 1.9475$ of the power law fitted to $L_p$
(Figure \ref{Lprimes} left) by least squares predicts that the counting function of $L_p$
will be asymptotic to a power law with exponent $1/b \approx 0.5135$. 
The power law fitted by least squares to the counting function of $L_p$ has exponent
$0.5115$ (Figure \ref{Lprimes} middle), different by only 0.002.

The variance function of $L_p$ is defined as 
the function $(0,\infty)\mapsto (0,\infty)$
from the mean of the first $n$ elements of $L_p$ 
to the variance of the first $n$ elements of $L_p$, for all $n\in\N$.
I estimated the variance function (Figure \ref{Lprimes} right) using
the first 25 elements of $L_p$, then the first 50, then the first 75,
and so on in successively longer intervals, each embedded in the next,
up to the first 2400 elements.
The exponent $1/b$ of an asymptotic regularly varying counting function
(Figure \ref{Lprimes} middle) predicts that the asymptotic variance function 
(cumulative) will be a power 
function with exponent 2 and coefficient $1/\bigl((1/b)(1/b+2)\bigr)= b^2/(1+2b)$.
The power law fitted by least squares to the variance function has exponent
approximately
$1.9831$ (Figure \ref{Lprimes} right), not greatly different from the predicted value 2.
The estimated coefficient is approximately $0.9676$, while $1/\bigl((1/b)(1/b+2)\bigr) \approx 0.7748$.

Since all primes are cyclics, if there are $k$ or more primes in $(n^2,(n+1)^2)$, 
then there are $k$ or more cyclics in $(n^2,(n+1)^2)$.

Analogous to $L_p$, define $L_c$ as
numbers $k$ such that the number of cyclics in every interval $[j^2,(j+1)^2], j>k$, 
is not less than the number of cyclics in the interval $[k^2,(k+1)^2]$.
Using the identical algorithm used to calculate $L_p$ for primes, with cyclics
replacing primes as the input, I calculated 769 values of $L_c$ based
on the 28,488,167 cyclics less than $10^8$.
For example, as the first 25 cyclics are  1, 2,    3,    5,    7,   11,   13,   15,   17,   
19,   23,   29,   31,   33,   35,   37,   41,   43,   47,   51,   53,   59,   61,   65,   67,
the seven intervals $[1, 4]$, $[4, 9],\ldots$, $[49, 64]$ contain  $N_c(n)= $ 2,    2,    3,    3,    4,    4,   and  4 cyclics,
and no later intervals \emph{in these calculations} have fewer than 4 cyclics. 
Consequently, I conjecture that $L_c(1) = 1, L_c(2) = 3, L_c(3) = 5$. At greater length, I conjecture:

\begin{conjecture}[$k$-fold Legendre for cyclics]\label{conj:Legendrek1cyclics}
For cyclics,
$L_c =(1, 3,          5,     8,       11,     14   , 15   , 16,
      19$,     $21,     27,     29$,    $ 33,     38,     39,     46,
      47,     51,     58,     61,     62,     66,     82,     86,
      90,    104,    105,    108,    110,    118,    126,    127,
     129$,    $131$,    $138,    141$,    $149,    152,    159,    161,
167,    170,    172,    174,    180,    182,    185,    187,\ldots ).$
\end{conjecture}

\begin{figure}[htb]
\centering
\includegraphics[width=1\textwidth]{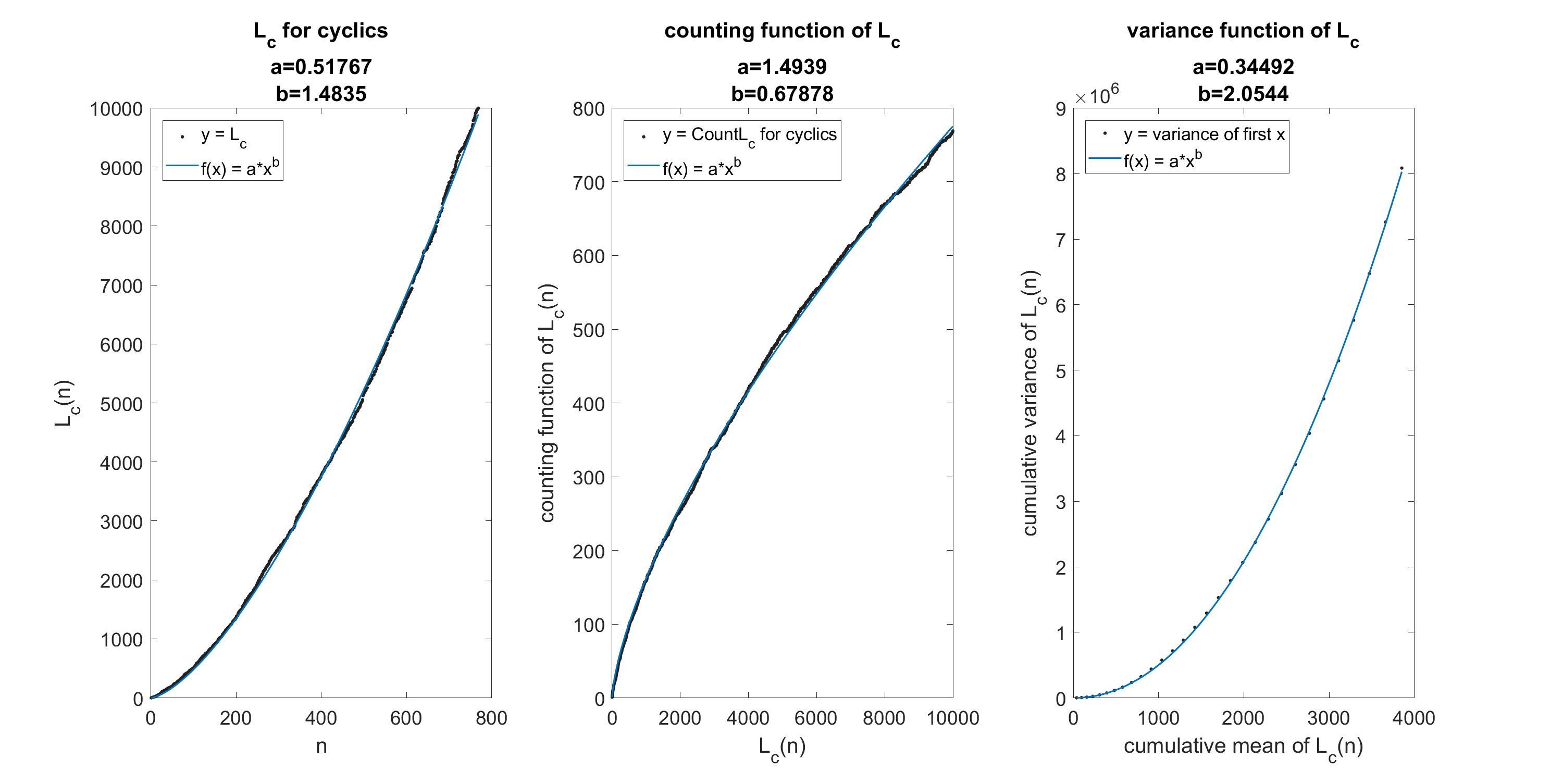}
\caption{(left) The sequence $L_c$ (black dots) 
and the power function (blue curve) fitted by least squares.
(middle) The counting function of $L_c$ (black dots) 
and the power function (blue curve) fitted by least squares.
(right) The variance function (black dots) of $L_c$ 
and a power function (blue curve)  fitted by least squares.}
\label{fig:Lcyclics}
\end{figure}

Figure \ref{fig:Lcyclics} (left) indicates that a power function approximates closely the
 calculated values of $L_c$ for cyclics. 
The counting function of $L_c$ (middle) 
and the variance function  of $L_c$ (right) 
are also well approximated by power functions.

\begin{conjecture}[asymptotic $k$-fold Legendre for cyclics]\label{conj:LegendrekasymptoticCyclics}
As $n\to\infty$, $L_c(n)$ is asymptotic to a regularly varying function 
with index that satisfies $1 \leq b  \leq 2$.
\end{conjecture}

\subsubsection{Landau's problem 4}

Fourth in Landau's list, the near-square conjecture states that infinitely many primes $p\in\Prime$ satisfy 
$p = n^2+1$ for some $n\in \N$.
These primes (\seqnum{A002496}) are one variety of the ``near-square primes.''
Hardy and Littlewood \cite{HardyLittlewood1923} conjectured a still unproved
asymptotic expression for the counting function of classes of primes including 
near-square primes.
Alfred Edward Western (1873--1961) \cite[p.\ 109, table]{Western1922} compared the conjectured asymptotic formula of Hardy and Littlewood with the numerically evaluated counting function of near-square primes
using a list \cite[pp.\ 238--239]{Cunningham1923}
by Allan J. C. Cunningham (1842--1928)
of the values of $n<15000$ such that $p = n^2+1$ is a near-square prime.
I programmed the computation and found, unexpectedly, that every one of Western's 15 tabulated values of the
counting function of near-square primes is one less than the corresponding value I obtained.
The discrepancy is due to Cunningham's omission \cite[p.\ 238]{Cunningham1923} of $n = 1$ for the first near-square prime, $1^2 + 1 = 2$.
Weisstein's list \cite{WeissteinNearSquare} of near-square primes $p = n^2+1$ begins with this case $n = 1$.

\begin{figure}[htb]
\centering
\includegraphics[width=1\textwidth]{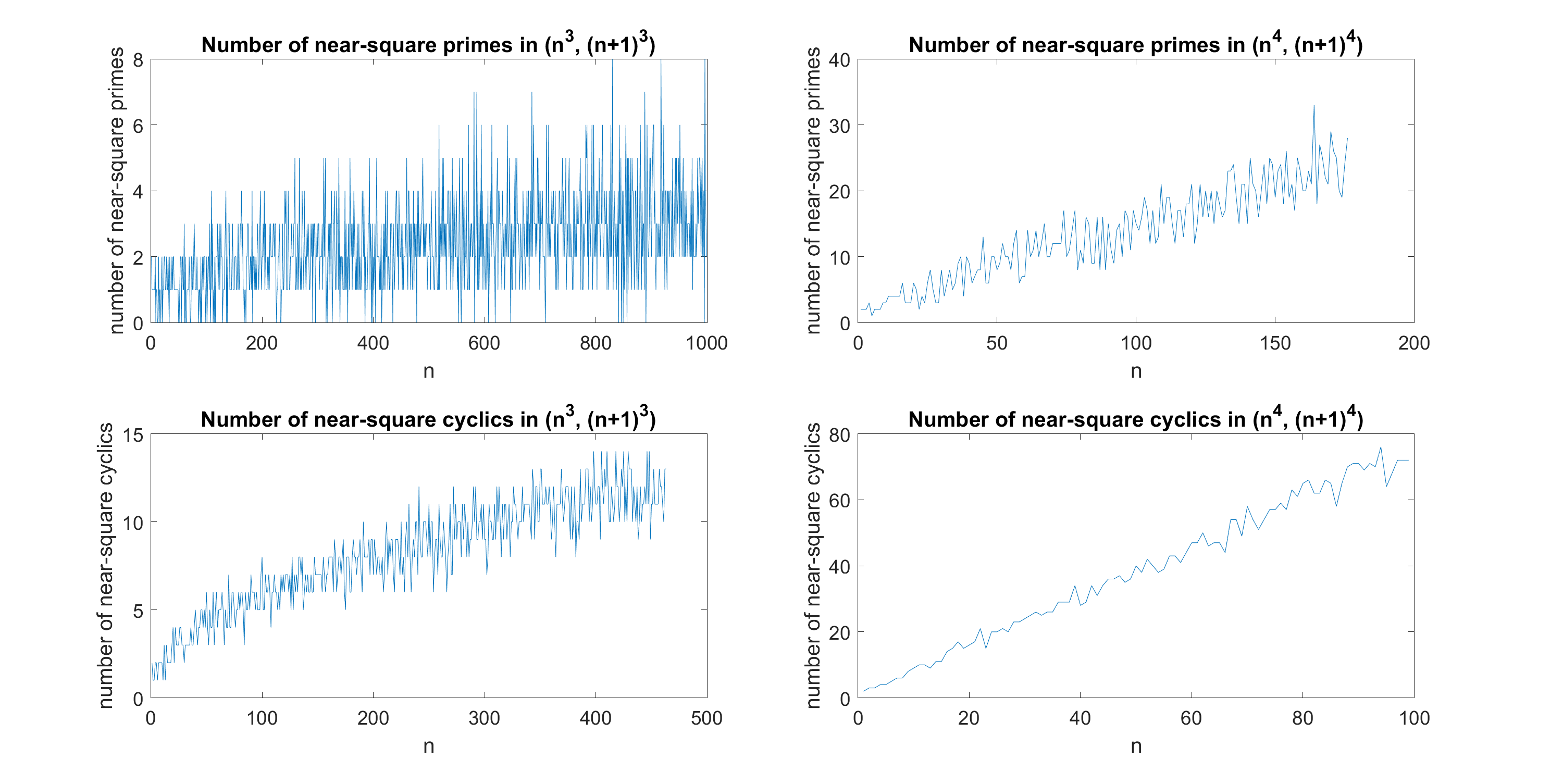}
\caption{(top row) The number of near-square primes $p=m^2+1$ in each interval (left) $(n^3, (n+1)^3)$ and (right) $(n^4, (n+1)^4)$.
(bottom row) The number of near-square cyclics $c=m^2+1$ in each interval (left) $(n^3, (n+1)^3)$ and (right) $(n^4, (n+1)^4)$.}
\label{fig:NearSquaresInIntervals3or4}
\end{figure}

W. A. Golubew (1891--1972) \cite[pp.\ 10--12]{Golubew1958} tabulated the values of $n\in [1, 10000]$ such that
$n^2+1\in\Prime$ (including $n=1$, unlike Cunningham)
and conjectured \cite[p.\ 13]{Golubew1958} that,
for every $m\in\N$, the interval $(m^4, (m+1)^4)$ contains at least one
near-square prime $p = n^2+1$ for some $n\in \N$.
I confirmed Golubew's conjecture using the 50,847,534 prime numbers less than $10^9$, 
which include 2379 near-square primes, the last being
$999444997 = 31614^2+1$ (Figure \ref{fig:NearSquaresInIntervals3or4}, top right).
Not every interval between successive cubes $(m^3, (m+1)^3)$ for $m\in\N$
contains at least one near-square prime (Figure \ref{fig:NearSquaresInIntervals3or4}, top left).
For example, the interval $(9^3 = 729,\ 10^3 = 1000)$ contains no near-square prime.

The 28,488,167 cyclics less than $10^8$ include 
3,786 near-square cyclics equal to $n^2+1$ for some $n\in \N$, beginning with
2, 5, 17, 37, 65, 101, 145, 197, 257, 401, 485, 577, 677, 785, 901, 1157, 1297,
1601, 1765, 1937,   2117, 2501,  2917,  3137,  3365,  3601,  3845,  4097,  4357,  
5477,  5777,  6085,  6401,  7057,  7397,  7745,  8101,  8465,  8837,  9217,  9605, 
10001, 10817, 11237, 11665, 12101, 12545, 12997, 13457, 14401, 14885, 15377, 
15877, 16385, 16901, 17957, 18497, 19045, 19601, 20165, 20737, 21317, 21905, 
22501, 23717, 24337, 24965, 25601, 26897, 27557, 28901, 30977, 31685, 32401, 
33857, 34597, 35345, 36101, 37637, 38417, 40001, 41617, 42437, 43265, 44101, 
45797, 46657, 48401, 49285, 50177, 51077, 52901, 54757, 55697, 56645, 57601, 
59537, 60517, 62501, 63505 
and ending with $99880037 = 9994^2+1$ and $99,920,017 = 9996^2+1$.

\begin{conjecture}[near-square analog for cyclics]\label{conj:near-square}
Infinitely many cyclics $c\in\C$ satisfy $c=n^2+1$ for some $n\in \N$.
\end{conjecture}

After seeing Conjecture \ref{conj:near-square} in an earlier draft, 
Carl Pomerance (personal communication, 2025-06-02 15:53) conjectured that
$n^2-1$ is cyclic for infinitely many $n$. He suggested that
the conjecture is plausible because all prime factors of $n^2-1$ are at most $n+1$.

\begin{conjecture}[near-square analog of Golubew for cyclics]\label{conj:near-squareGolubew}
For every $m\in\N$, the interval $(m^3, (m+1)^3)$ contains at least one
near-square cyclic $c = n^2+1$ for some $n\in \N$.
For every $m\in\N$, the interval $(m^4, (m+1)^4)$ contains at least two
near-square cyclics $c = n^2+1,\ c' = n'^2+1$ for some $n<n'\in \N$.
\end{conjecture}

\subsection{Oppermann's conjecture}
Ludvig Henrik Ferdinand Oppermann (1817--1883) \cite[p.\ 174]{oppermann1882}, in an unpublished lecture on March 9, 1877, conjectured
that for every $n>1$, there exist two primes $p,\ p'$
such that $n^2 - n < p < n^2 < p' <n^2 + n$. 
Several claims to prove Oppermann's conjecture have been published or posted
but I am not aware that any has been independently confirmed. 

The heuristic approach of Deligne's Conjecture \ref{conj:Np} suggests that,
asymptotically as $n\to\infty$, each number in the interval $[n^2 - n ,n^2]$
has probability $1/\log(n^2)$ of being prime, 
hence the number of primes in $[n^2 - n ,n^2]$ should be asymptotic to
$n/(2\log(n))$.
The same argument and conjectured conclusion hold for the number of primes in $[n^2  ,n^2+n]$.
These suggestions make plausible this conjecture:

\begin{conjecture}[counting Oppermann primes]\label{conj:Opp1}
As $n\to\infty$, the number
of primes in $[n^2 - n ,n^2]$ is asymptotic to $n/(2\log n)$
and the number of primes in $[n^2  ,n^2+n]$ is also asymptotic to $n/(2\log n)$. 
\end{conjecture}

If true, Conjecture \ref{conj:Opp1} would imply:
\begin{conjecture}[$k$-fold Oppermann conjecture for primes]\label{conj:Oppermannkfoldprimes}
For every $k\in\N$, there exists $N(k)\in\N$ such that for all $n > N(k)$, 
there exist at least $k$ primes in $[n^2 - n, n^2]$ and another at least 
$k$ primes in $[n^2, n^2 + n]$.
In particular, $N(2) = 16$, $N(3) =     36$, $N(4) =     46$, $N(5) =      76$, 
$N(6) =  N(7) =  79$,   $N(8) = 85$, $N(9) = 118$,      $N(10) =  136$,  
$ N(11) = N(12) =     155$, $N(13) = 188$.
\end{conjecture}
Oppermann's conjecture corresponds to $N(1) = 1$.

\begin{remark}
Oppermann's conjecture for primes implies Legendre's conjecture for primes.
Specifically, if (a) for every $n\in\N,\ n>1$, there exists $p\in\Prime$
such that $n(n-1)<p<n^2$, or 
(b) for every $n\in\N$, there exists $p'\in\Prime$
such that $n^2<p'<n(n+1)$, or (c) both (a) and (b) hold, then
(d) for every $n\in \N$, there exists $p\in\Prime$
such that $n^2<p<(n+1)^2$. 
\end{remark}

\begin{proof}
Since $(n-1)^2 < n(n-1)$ for all $n>1$, (a) implies (d).
Since $n(n+1) < (n+1)^2$ for all $n\in\N$, (b) implies (d).
\end{proof}

\begin{conjecture}[Oppermann analog for cyclics]\label{conj:Oppermann}
For every $n > 1$, there exist $c,\ c'\in\C$
such that $n^2 - n < c < n^2 < c' <n^2 + n$. 
The number of cyclics in $[n^2 - n ,n^2]$ is asymptotic to
\begin{equation}\label{eq:CyclicsInSquaresUsingPollackLeft}
 \frac{n}{e^{\gamma}\log \log \log (n^2)}\left(1 - \frac{{\gamma}}{\log \log \log (n^2)} \right)
\sim  \frac{n}{e^{\gamma}\log \log \log (n)}\left(1 - \frac{{\gamma}}{\log \log \log (n)} \right).
\end{equation}
The number of cyclics in the interval $[n^2  ,n^2+n]$ is 
asymptotic to 
\begin{equation}\label{eq:CyclicsInSquaresUsingPollackRight}
 \frac{n}{e^{\gamma}\log \log \log (n(n+1))}\left(1 - \frac{{\gamma}}{\log \log \log (n(n+1))} \right)
\sim  \frac{n}{e^{\gamma}\log \log \log (n)}\left(1 - \frac{{\gamma}}{\log \log \log (n)} \right).
\end{equation}
\end{conjecture}

The sum of \eqref{eq:CyclicsInSquaresUsingPollackLeft} and
\eqref{eq:CyclicsInSquaresUsingPollackRight} is asymptotic to 
\eqref{eq:CyclicsInSquaresUsingPollack}, and 
\eqref{eq:CyclicsInSquaresUsingPollackLeft}
is asymptotic to 
\eqref{eq:CyclicsInSquaresUsingPollackRight}.
For $n\in\N$ such that $n\geq 4$, one has 
\eqref{eq:CyclicsInSquaresUsingPollackLeft}
$<$ \eqref{eq:CyclicsInSquaresUsingPollackRight}.

\begin{figure}[htb]
\centering
\includegraphics[width=1\textwidth]{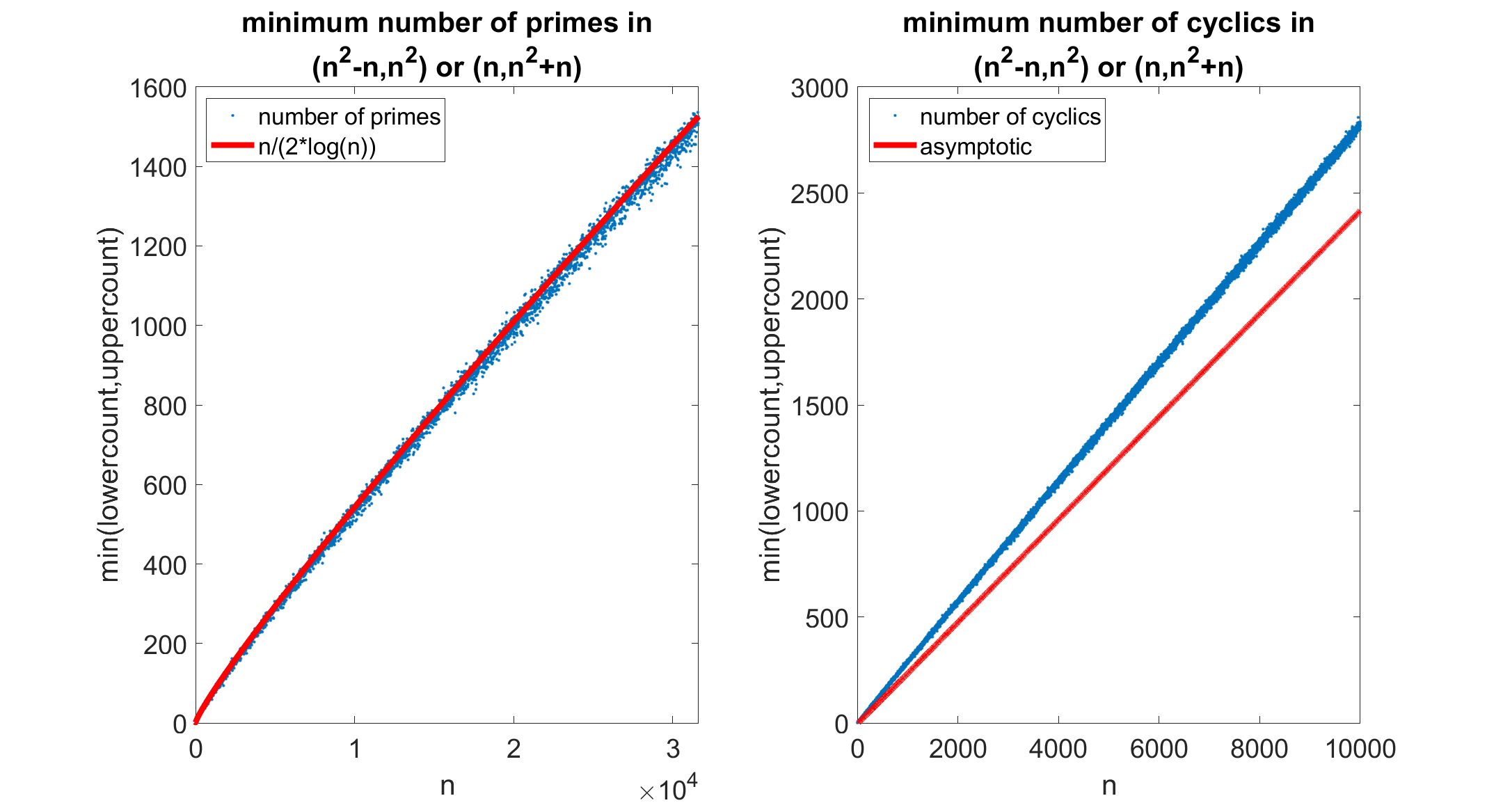}
\caption{(left) For the primes less than $10^9$, 
the minimum of the numbers of primes in the two intervals
$[n^2 - n, n^2]$ and $[n^2, n^2 + n]$ is shown by a blue dot
for every tenth value of $n$ 
to avoid having the blue dots overwrite the red curve.
The red curve plots $n/(2 \log n),\ n = 2,\ldots,31622$.
(right) For the cyclics less than $10^8$, for each $n = 20,\ 21,\ 22,\ldots,9999$, a blue dot shows
the minimum of the numbers of cyclics in the two intervals
$[n^2 - n, n^2]$ and $[n^2, n^2 + n]$.
The red line shows the asymptotic expression in 
\eqref{eq:CyclicsInSquaresUsingPollackLeft}.}
\label{fig:OppermannPrimesCyclics}
\end{figure}

For the primes less than $10^9$, for each $n = 2,\ldots,31622$, Figure \ref{fig:OppermannPrimesCyclics}  (left)
plots $n/(2 \log n)$ and
the minimum of the numbers of primes in the two intervals
$[n^2 - n, n^2]$ and $[n^2, n^2 + n]$.
For the cyclics less than $10^8$,
Figure \ref{fig:OppermannPrimesCyclics}  (right) 
compares the asymptotic expression in 
\eqref{eq:CyclicsInSquaresUsingPollackLeft} (red curve) with 
the minimum of the numbers of cyclics in the two intervals
$[n^2 - n, n^2]$ and $[n^2, n^2 + n]$ (blue dot).

If true, Conjecture \ref{conj:Oppermann} would imply:
\begin{conjecture}[$k$-fold Oppermann conjecture for cyclics]\label{conj:Oppermannkfoldcyclics}
For every $k\in\N$, there exists $N(k)\in\N$ such that for all $n > N(k)$, 
there exist at least $k$ cyclic numbers in the interval $[n^2 - n, n^2]$ and another at least 
$k$ prime numbers in the interval $[n^2, n^2 + n]$.
In particular, $N(2) = 4$, $N(3) =     7$, $N(4) =     13$, $N(5) =      16$, 
$N(6) = 18$, $N(7) =  21$,   $N(8) = 25$, $N(9) = 31$,      $N(10) = N(11) = 32$,  
$  N(12) =     40$, $N(13) = 44$.
\end{conjecture}

 Numerically, there exist $c, c'\in\C$ 
such that $n^2 - n < c < n^2 < c' <n^2 + n$  for each $1< n\leq 9998$.
For $n=2, 3,\ldots,25$, the lesser of the number of cyclics in $(n^2 - n, n^2)$ and 
the number of cyclics in $(n^2, n^2 + n)$ is 
1,    1,    2,    1,    2,    2,    2,    3,    3,    3,    3,    4,    3,    5,    5,    4,    6,    5,    6,    7,    6,    7,    9,    7,   and  7.
All computed later elements in this sequence are 8 or larger.

\subsection{Brocard's and Desboves' conjectures}

Pierre Ren\'e Jean Baptiste Henri Brocard (1845--1922) \cite{brocard1904} conjectured in 1904 that
there are at least four primes between the squares of two successive primes,
provided that the first prime be greater than 3.
However, between the squares 9 and 25 of 3 and 5, respectively, there are more
than four primes, namely, 11,    13,   17,    19,    23,
so Brocard's proviso should be replaced by requiring that the first prime
be greater than 2.
Several claims to prove Brocard's  conjecture have been published or posted
but I am not aware that any has been independently confirmed. 

\begin{remark}
Desboves' conjecture \cite{Desboves1855}
that for every $n\in \N$, there exist $p, p'\in\Prime$ 
such that $n^2<p < p'<(n+1)^2$ implies
Brocard's (adjusted) conjecture \cite{brocard1904} that
there are at least four primes between the squares of two successive primes greater than 2.
More generally, if $p,\ p'$ are two primes both greater than 2,
then Desboves' conjecture \cite{Desboves1855} implies that there are at least
$2(|p'-p|)$ primes between $p^2$ and $p'^2$.
\end{remark}

\begin{proof}
As 2 is the only even prime, every pair of consecutive primes except 2 and 3 is separated by at least one even number. 
If $p>2$ and $p'=p+2$ are twin primes and $p<m<p'$, $m\in\N$, then, 
by Desboves' conjecture, there are at least two primes between $p^2$ and $m^2$.
Again by Desboves' conjecture, there are at least two primes between $m^2$ and $p'^2$.
So there are at least four primes between $p^2$ and $p'^2$.
All other pairs $2<p < p'$ of consecutive primes have more than one intervening 
even number and one or more intervening odd numbers between them, 
and each of those intervening numbers contributes two or more primes to the
count of primes between $p^2$ and $p'^2$.
\end{proof}

For $k\in\N,\ k\geq 4$, define
$B(k)$ to be the smallest $n\in\N$ such that,
for all $m\in\N$ with $m\geq n$, there are always at least 
$k$ primes between $p_m^2$ and $p_{m+1}^2$.
For example, Brocard's conjecture asserts that $B(4)=2$.
(This 2 points to the second prime, $p_2 = 3$, not to the first prime 2.)
In \seqnum{A050216}, T. D. Noe calculated
the numbers of primes between $p_n^2$ and $p_{n+1}^2$  
for $n\leq 10^4$.
Based on Noe's calculations, I conjecture:

\begin{conjecture}[$k$-fold Brocard for primes]\label{conj:Brocardk}
Let $B(4) =2, B(5)=2, B(6) = 3, B(7) =B(8) = B(9)= 5, B(10) =B(11) = 7, B(12)=B(13) =
B(14)=B(15)=B(16) = 10, B(17) = B(18) = B(19) = B(20) = 13$.
Then for $k=4,\ldots, 20$ and for every $n\in \N$ such that
$n \geq B(k)$, there exist at least $k$ primes 
between $p_n^2$ and $p_{n+1}^2$.
More generally, for every $k\in\N$, there exists $B(k)\in\N$
such that for every $n\in \N$ with $n \geq B(k)$, there exist at least $k$ 
primes between $n^2$ and $(n+1)^2$.
\end{conjecture}

\begin{conjecture}[Brocard analog for cyclics]\label{conj:Brocard}
For every $n > 2$, there exist at least six cyclics $c\in\C$
such that $c_n^2 < c < c_{n+1}^2$.
\end{conjecture}
 There are two cyclics $c_2 = 2 ,\ c_3 = 3$ 
in the open interval $(c_1^2 = 1, c_2^2 = 4)$.
There are two cyclics $c_4 = 5$, $c_5 = 7$
in $(c_2^2 = 4, c_3^2 = 9)$.
The six cyclics in $(c_{3}^2 = 9, c_4^2 = 25)$ are
 11,    13,    15,    17,    19, and    23.
For $n = 1,\ldots, 25$, the number of cyclics $c\in\C$
such that $c_n^2< c < c_{n+1}^2$ is
2,    2,    6,     8,    25,    16,    17,    21,    22,    56,   102,    36,    36,
45,    49,    96, 52,   113,   125,    65,   206,    80,   152,    83,    84.
For $3 \leq n\leq 1009$ (with $c_{1009} = 3157$),
there are at least six cyclics $c\in\C$
such that $c_n^2 < c < c_{n+1}^2$.

If the $k$-fold Brocard conjecture for primes is true, then it is immediate
that for $k=4,\ldots, 20$ and for every $n\in \N$ such that $n \geq B(k)$, there exist at least $k$ cyclics 
between $p_n^2$ and $p_{n+1}^2$, since all primes are cyclics.
The following analog for cyclic numbers is not an obvious consequence of 
the $k$-fold Brocard conjecture for primes,
and may hold true even if Conjecture \ref{conj:Brocardk} is false.

\begin{conjecture}[$k$-fold Brocard analog for cyclics]\label{conj:Brocardkcyclic}
For every $k\in\N$, there exists $C(k)\in\N$ such that for every $n\in \N$
with $n \geq C(k)$, there exist at least $k$ cyclic numbers $c$
such that $c_n^2 < c < c_{n+1}^2$.
\end{conjecture}
Numerically, $C(2) = 1$, i.e., for every $n\geq 1$ computed here,
there exist at least $2$ cyclic numbers $c$
such that $c_n^2 < c < c_{n+1}^2$.

\subsection{Schinzel's conjectures}

Andrzej Bobola Maria Schinzel (1937--2021) \cite[p.\ 155, Conjecture P$_1$]{Sierpinski1988} conjectured in 1961 that,
for real number $x\geq 117$, there is at least one prime between $x$ and
$x+\sqrt{x}$. 
When $x$ is evaluated only at primes, the exceptional cases
(among the primes less than $10^9$, in my computations)
where there is \emph{no} prime between $p_n$ and  $p_n+\sqrt{p_n}$ are
$p_2 = 3, p_4 = 7, p_6 = 13, p_9 = 23, p_{11} = 31$, and $p_{30} = 113$.
A stronger conjecture by Schinzel \cite[p.\ 156]{Sierpinski1988} is that, 
for real number $x\geq 8$, there is at least one prime between 
$x$ and $x+(\log{x})^2$.
When $x$ is evaluated only at primes, the exceptional cases 
(among the primes less than $10^9$) where
there is \emph{no} prime between $p_n$ and  $p_n+(\log{p_n})^2$ are
$p_1 = 2, p_2 = 3$, and $p_4 = 7$.
With their different lower bounds, both conjectures have been confirmed
numerically for $x \leq 4.44 \times 10^{12}$.
Conjecture P$_1$ implies Conjecture P  of Wac\l{}aw Sierpi\'nski (1882--1969)
 \cite[p.\ 153]{Sierpinski1988}. 
Both of Schinzel's conjectures have obvious analogs for cyclics
and one analog not so obvious.
I verified these conjectures for the cyclics less than $10^8$.

\begin{conjecture}[Schinzel Conjecture P$_1$ analog for cyclics]\label{conj:SchinzelP1}
For every $n\in\N$, 
$$ c_{n+1} \leq c_n+ \sqrt{c_n} $$
except for $c_3 = 3, c_5 = 7$, and $c_{11} = 23$.
\end{conjecture}

 \begin{conjecture}[Schinzel conjecture for log$^2$ analog for cyclics]\label{conj:Schinzellog2}
For every $n\in\N$, 
$$ c_{n+1} \leq c_n+ (\log{c_n})^2 $$
except for $c_1 = 1, c_2 = 2, c_3 = 3$, and $c_{5} = 7$.
\end{conjecture}

\begin{conjecture}[Schinzel-type conjecture for $2\times\log$ analog for cyclics]\label{conj:Schinzel2log}
For every $n\in\N$, 
$$ c_{n+1} \leq c_n+ 2\log{c_n} $$
except for $c_1 = 1$ and $c_{5} = 7$.
\end{conjecture}

\subsection{Golubew's conjectures}

Noting Legendre's conjecture for primes,
Golubew \cite[p.\ 85]{Golubew1957} conjectured in 1957 that for $n\in\N$, 
there is at least one pair of twin primes $p,\ p+2$ between $n^3$ and $(n +1)^3$.
Further, he conjectured that, for $n\in\N$, 
there is at least one quartet of primes $p,\ p+2,\ p+6,\ p+8$ between $n^5$ and $(n +1)^5$.
For $n=1,\ldots, 25$, I find that the number of pairs of twin primes between $n^3$ and $(n +1)^3$ is
2,  2,   3,   3,   5,   5,   4,   6,   5,  11,   9,  12,  11,  12,  17,  17,  16,  19,  16,  18,  24,  22,  17,  22,  and 26.
For $n=1,\ldots, 10^3-1$ such that $(n +1)^3\leq 10^9$,
I find numerically that the number of pairs of twin primes between $n^3$ and $(n +1)^3$ is never less than two.
As examples, if $n=1$, then $(3,5)$ and $(5,7)$ are two pairs of twin primes between 1 and $(n+1)^3=8$;
and if $n=2$, then $(11,13)$ and $(17,19)$ are two pairs of twin primes between 8 and 27.
Moreover, I find 3 or more pairs of twin primes  between $n^3$ and $(n +1)^3$ for all $3\leq n \leq 999$,
4 or more pairs of twin primes  between $n^3$ and $(n +1)^3$ for all $5\leq n \leq 999$, 
5 or more pairs of twin primes  between $n^3$ and $(n +1)^3$ for all $8\leq n \leq 999$, 
6 or more pairs of twin primes  between $n^3$ and $(n +1)^3$ for all $10\leq n \leq 999$, and so on.

Golubew \cite[p.\ 84, Table 2]{Golubew1957} tabulated the number of pairs of twin primes between $n^3$ and $(n +1)^3$ for $n=1,\ldots, 80$.
I confirmed his counts with five exceptions, which I believe to be his errors.
I list the five cases in which I challenge his results with three numbers:
 $n$, his count of twin primes between $n^3$ and $(n +1)^3$, and my count of twin primes between $n^3$ and $(n +1)^3$.
These five cases are:
25, 27, 26;
26, 31, 32;
70, 109, 119;
74, 130, 131;
80, 160, 161.

\begin{conjecture}[number of twin primes between consecutive cubes]\label{conj:twinprimesbetwcubes}
For every $n\in\N$, 
the number of pairs of twin primes between $n^3$ and $(n +1)^3$ is never less than two.
More generally, for every $k\in\N$, there exists $N(k)\in\N$ such that for all $n\geq N(k)$ there are at least $k$ pairs of twin primes between $n^3$ and $(n +1)^3$.
Specifically, $N(1)=N(2) =1$, $N(3)=3$, $N(4) = 5$, $N(5)=8$,
$N(6)=N(7)=N(8)=N(9) = 10$, $N(10)= 11$, $N(11)=13$, $N(12) = N(13)=N(14)=N(15)=N(16) = 15$, and
$N(17) = 20$.
\end{conjecture}

Two consecutive primes $p, p'$ are defined \cite[p.\ 336]{Wolf1998} to be cousin primes if $|p-p'| = 4$ and defined to be sexy primes if $|p-p'| = 6$.
By these definitions, 3 and 7 are not cousin primes and 11 and 17 are not sexy primes because they are not consecutive. 
Every other pair of primes $p, p'$ with $|p-p'| = 4$ is consecutive, hence cousin.
Many other pairs of primes $p, p'$ with $|p-p'| = 6$ are not consecutive, hence not sexy.
For $n=1,\ldots, 25$, I find that the number of pairs of cousin primes between $n^3$ and $(n +1)^3$ is
    0,   2,   2,   5,   3,   5,   8,   3,  11,   7,  12,   7,  15,  14,  13,  10,  19,  13,  20,  21,  22,  23,  24,  28, and 31.
For $n=1,\ldots, 25$, I find that the number of pairs of sexy primes between $n^3$ and $(n +1)^3$ is
   0,   0,   3,   2,   5,   6,   7,  11,   7,  15,  11,  12,  19,  15,  20,  21,  30,  27,  29,  33,  30,  37,  43,  36,  and 52.
These numbers and all the following counts of cousin primes between successive cubes suggest:

\begin{conjecture}[number of cousin primes between consecutive cubes]\label{conj:cousinprimesbetwcubes}
For  $n\in\N$ with $n>1$, 
the number of pairs of cousin primes between $n^3$ and $(n +1)^3$ is never less than two.
More generally, for every $k\in\N$, there exists $N(k)\in\N$ such that,
for all $n\geq N(k)$, 
there are at least $k$ pairs of cousin primes between $n^3$ and $(n +1)^3$.
Specifically, $N(1)=N(2) =2$, $N(3)=8$, $N(4) = N(5) = N(6) = N(7) = 9$, $N(8)=N(9)=N(10)=12$.
\end{conjecture}

\begin{conjecture}[number of sexy primes between consecutive cubes]\label{conj:sexyprimesbetwcubes}
For  $n\in\N$ with $n>2$, 
the number of pairs of sexy primes between $n^3$ and $(n +1)^3$ is never less than two.
(While 3 pairs of sexy primes occur between $3^3$ and $4^3$, only 2 pairs of sexy primes occur between $4^3$ and $5^3$.)
More generally, for every $k\in\N$, there exists $N(k)\in\N$ such that,
for all $n\geq N(k)$, there are at least $k$ pairs of sexy primes between $n^3$ and $(n +1)^3$.
Specifically, $N(1)=N(2) =3$, $N(3)=N(4) = N(5) = 5$, $N(6) =6$, $N(7) = 7$, $N(8)=N(9)=N(10)=N(11) = 11$.
\end{conjecture}

\begin{conjecture}[asymptotic $k$-fold primes between consecutive cubes]\label{conj:Golubewkasymptoticprimes}
As $n\to\infty$, the numbers of pairs of twin primes, cousin primes, and sexy primes between consecutive cubes $n^3$ and $(n +1)^3$
are asymptotic to regularly varying functions of $n$
with indices between $3/2$ and $2$, and the indices for twin primes and cousin primes are identical.
\end{conjecture}

\begin{figure}[htb]
\centering
\includegraphics[width=1\textwidth]{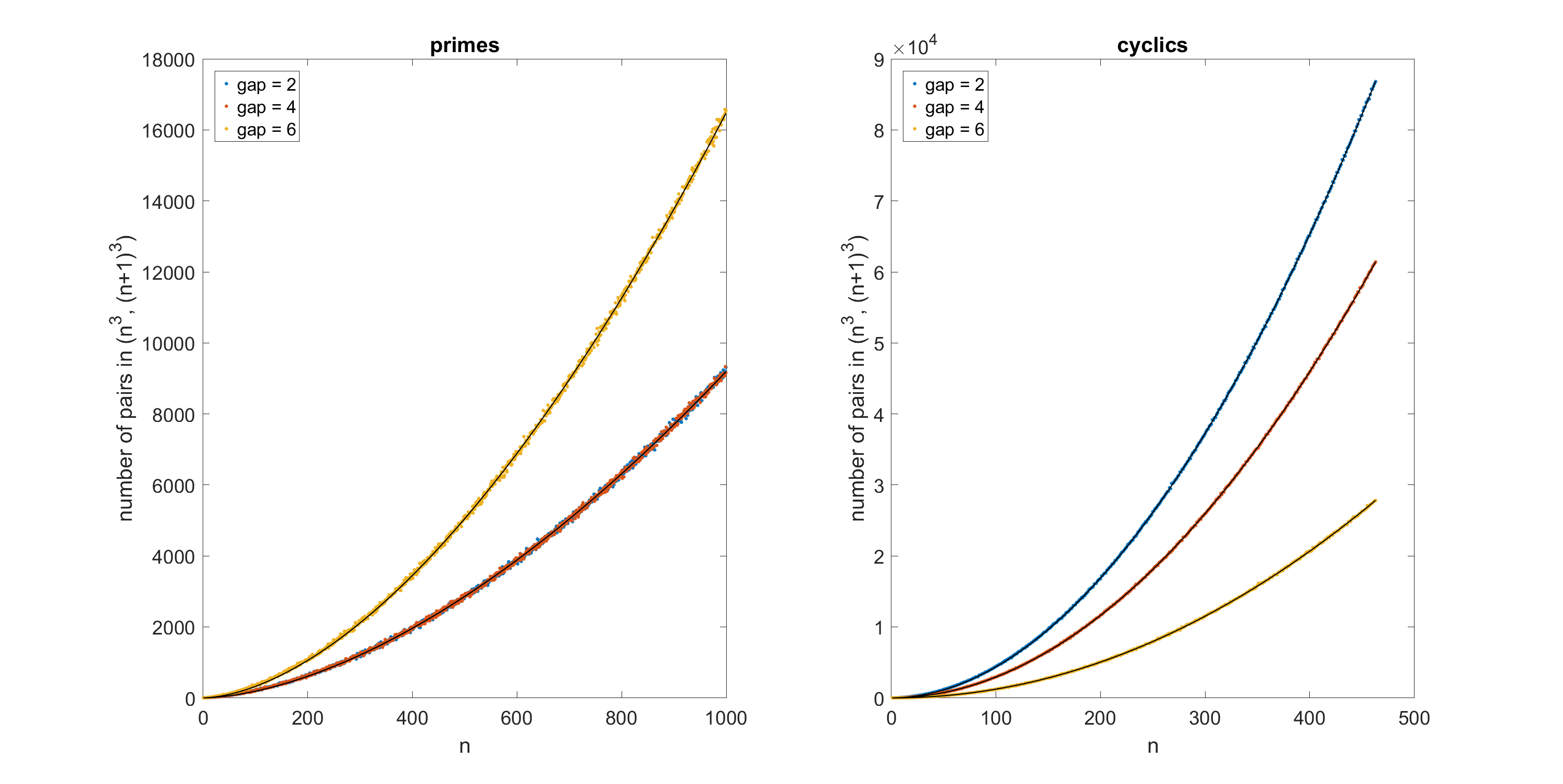}
\caption{(left) Each dot represents, on the vertical axis, the number of pairs of consecutive primes $p_m, p_{m+1}$ in the interval between $n^3$ and $(n +1)^3$
for the value of $n$ on the horizontal axis.
Gaps $p_{m+1}-p_m=2, 4, 6 $ correspond to twin (blue dots), cousin (red dots), and sexy (yellow dots) primes.
Black curves show power functions $an^{b}$ fitted by least squares.
For twin primes (gap = 2), $a = 0.079, b = 1.689$. 
For cousin primes (gap = 4), $a = 0.0804, b = 1.686$.
For sexy primes (gap = 6), $a = 0.124, b = 1.708$.
(right)  Each dot represents, on the vertical axis, the number of pairs of consecutive cyclics $c_m, c_{m+1}$ in the interval between $n^3$ and $(n +1)^3$
for the value of $n$ on the horizontal axis.
Gaps $c_{m+1}-c_m=2, 4, 6 $ correspond to twin (blue dots), cousin (red dots), and sexy (yellow dots) cyclics.
Black curves show power functions $an^{b}$ fitted by least squares.
For twin cyclics (gap = 2), $a = 0.5493, b = 1.950$. 
For cousin cyclics (gap = 4), $a = 0.3187, b = 1.983$.
For sexy cyclics (gap = 6), $a = 0.1051, b = 2.034$.}
\label{fig:GolubewPrimesCyclics}
\end{figure}

Figure \ref{fig:GolubewPrimesCyclics} (left) supports Conjecture \ref{conj:Golubewkasymptoticprimes} numerically.

Turning from primes to cyclics: twin cyclics, cousin cyclics, and sexy cyclics are pairs of consecutive cyclics with first difference 2, 4, and 6, respectively.
For $n=1,\ldots, 25$, I find that the numbers of pairs of twin cyclics between $n^3$ and $(n +1)^3$ are
         2,      4,       7,    13,    17,     22,    32,    41,    44,    57,    70,    80,    99,   107,   122,    132,   142,   171,   189,   220,   221,   239,   271,   292,   and 310.
For $n=1,\ldots, 25$, I find that the numbers of pairs of cousin cyclics between $n^3$ and $(n +1)^3$ are
0,   1,   3,   8,   8,   14,   15,   22,   29,   37,   36,   51,   50,   69,   71,   95,   92,   97,   120,   129,   142,   149,   177,   175,   and 194.
For $n=1,\ldots, 25$, I find that the numbers of pairs of sexy cyclics between $n^3$ and $(n +1)^3$ are
0,   0,   1,   0,   2,   4,   7,   7,   9,   8,   13,   13,   19,   17,   16,   23,   38,   44,   36,   42,   46,   58,   54,   67,  and 70.
These numbers and all the following counts of cyclics between successive cubes $(n^3, (n+1)^3)$ suggest:

\begin{conjecture}[number of twin cyclics between consecutive cubes]\label{conj:twincyclicsbetwcubes}
For every $n\in\N$, 
the number of pairs of twin cyclics between $n^3$ and $(n +1)^3$ is never less than two.
For every $k\in\N$, there exists $N(k)\in\N$ such that for all $n\geq N(k)$ there are at least $k$ pairs of twin cyclics between $n^3$ and $(n +1)^3$.
Specifically, $N(1)=N(2) =1$, $N(3)=N(4) = 2$, $N(4) = N(5) = N(6) = N(7) = 3$, $N(8) =N(9)=N(10)=N(11)=N(12)=N(13) = 4$.
\end{conjecture}

\begin{conjecture}[number of cousin cyclics between consecutive cubes]\label{conj:cousincyclicsbetwcubes}
For  $n\in\N$ with $n>1$, 
the number of pairs of cousin cyclics between $n^3$ and $(n +1)^3$ is never less than one.
For every $k\in\N$, there exists $N(k)\in\N$ such that for all $n\geq N(k)$ there are at least $k$ pairs of cousin cyclics between $n^3$ and $(n +1)^3$.
Specifically, $N(1)=2$, $N(2) =N(3)=3$, $N(4) = N(5) = N(6) = N(7) =N(8) = 4$, $N(9)=N(10)=N(11)=N(12)=N(13)=N(14)=6$.
\end{conjecture}

\begin{conjecture}[number of sexy cyclics between consecutive cubes]\label{conj:sexycyclicsbetwcubes}
For  $n\in\N$ with $n\geq 5$, 
the number of pairs of sexy cyclics between $n^3$ and $(n +1)^3$ is never less than two.
For every $k\in\N$, there exists $N(k)\in\N$ such that for all $n\geq N(k)$ there are at least $k$ pairs of sexy cyclics between $n^3$ and $(n +1)^3$.
Specifically, $N(1)=N(2) =5$, $N(3)=N(4) =  6$, $N(5) =N(6) =N(7) = 7$, $N(8)=10$, $N(9)=N(10)=N(11) = N(12)=N(13) = 11$.
\end{conjecture}

\begin{conjecture}[asymptotic $k$-fold cyclics between consecutive cubes]\label{conj:Golubewkasymptoticcyclics}
As $n\to\infty$, the numbers of twin cyclics, cousin cyclics, and sexy cyclics between consecutive cubes $n^3$ and $(n +1)^3$
are asymptotic to regularly varying functions of $n$
with indices between 1 and $5/2$.
\end{conjecture}

\subsection{Sophie Germain primes and cyclics}
In approaching a proof of Fermat's Last Theorem,
Marie-Sophie Germain (1776--1831) considered
(as a very special case of much more general hypotheses) 
pairs of primes $(p,2p+1)$ such as $(3,7)$ and $(5,11)$.
A prime $p$ such that $2p+1$ is prime is called a Sophie Germain prime
(or an SG prime), and the prime $2p+1$ is called a safe or auxiliary prime.
Identifying the first appearance historically of SG primes is challenging
because much of Germain's work was never published, appeared 
in correspondence, is mentioned in the work of other mathematicians, 
or was published anonymously or pseudonymously
\cite{LaubenbacherPengelley2010}.

It is conjectured but unproved that there are infinitely many SG primes.
For positive real $x$, let $\pi_{SG}(x)$ be the counting function of SG primes,
i.e., the number of SG primes less than or equal to $x$.
It is conjectured but not proved \cite[p.\ 123]{Shoup2009} that, as $x\to\infty$,
\begin{equation*}
\pi_{SG}(x) \sim \left(2\prod_{\{p\in\Prime \mid p>2\}}\frac{p(p-2)}{(p-1)^2} \right)\frac{x}{(\log x)^2}.
\end{equation*}

If there are infinitely many SG primes, then Conjecture 
\ref{conj:GermainCyclics} is true, 
but Conjecture \ref{conj:GermainCyclics} may be true even if 
there are only finitely many SG primes.

Define a cyclic $c\in\C$ to be a Sophie Germain cyclic (or an SG cyclic) if
$2c+1\in\C$.

\begin{conjecture}[infinitely many SG cyclics]\label{conj:GermainCyclics}
There are infinitely many Sophie Germain cyclics.
\end{conjecture}

After seeing Conjecture \ref{conj:GermainCyclics} in a draft of this paper, 
Carl Pomerance (personal communication, 2025-05-24 11:02) announced
that he can prove that the number of Sophie Germain cyclics less than $x$
is asymptotic to $cx(e^{\gamma}\log \log \log x)^{-2}$ for an
appropriate $c>0$, and hence that Conjecture \ref{conj:GermainCyclics} is true.

Let $\sigma_n$ be the $n$th SG cyclic. The first 25 SG cyclics are
   1,    2,  3,  5,  7,   11,   15,   17,   23,   29,   33,   35,   41,   43,  
 47,   51,   53,   59,   61,   65,   69,   71,   79,   83,   89.
For example, $\sigma_7 = 15$ is an SG cyclic because $2\times 15 + 1=31\in \Prime$ is cyclic, although 31 itself is not an SG cyclic, 
as the following list shows.
The first 25 cyclics that are {\em not} SG cyclics are
13,    19,   31,   37,   67,   73,   77,   85,   87,   91,   97,
   101,   103,   109,   115,   137,   139,   145,   157,   163,   177,   181,
187,   193,   199. For example, $2\times 13 + 1=27$.

It is well known that every SG prime except SG prime 3 is congruent to 
2 mod 3, for if an SG prime $p$ were congruent to 1 mod 3, then
$2p+1$ would be congruent to 3 mod 3, i.e., composite, 
contradicting the assumption that $p$ is an SG prime.
The SG cyclics are different.

\begin{conjecture}[SG cyclics mod 3]\label{conj:SGcyclicsMod3}
As the number of SG cyclics grows without limit,
the number of SG cyclics congruent to $j$ mod 3, $j=1, 2, 3$, grows without limit, and the limiting fraction of SG cyclics congruent to 1 mod 3 equals
the limiting fraction of SG cyclics congruent to 3 mod 3.
\end{conjecture}

After seeing an earlier version of Conjecture \ref{conj:SGcyclicsMod3}, 
Carl Pomerance (personal communication, 2025-05-24 11:02) announced
that he can show that the fractions of SG cyclics congruent to $j$ mod $3$ 
approach the limits $w_1=w_3=0$ and $w_2=1$.

Based on the first 3,441,316 SG cyclics,
the proportions of SG cyclics congruent to 1, 2, and 3 mod 3 are 
approximately 0.1360, 0.7252, 0.1388.
Based on the first 6,882,632 SG cyclics (twice as many),
the proportions of SG cyclics congruent to 1, 2, and 3 mod 3 are 
approximately 0.1342, 0.7290, 0.1368.

\begin{conjecture}[Desboves analog for SG cyclics]\label{conj:LegendreSGcyclics}
For every $n\in \N$, there exists at least two SG cyclics in $(n^2, (n+1)^2)$.
\end{conjecture}

For $n=1,\ldots,7070$, every interval $(n^2,(n+1)^2)$ contains
at least two SG cyclics. For example, for $n=1,\ldots,25$,
the number of SG cyclics in each interval $(n^2,(n+1)^2)$ is
2,   2,   2,   2,   3,   3,   4,   4,   3,   3,   6,   5,   4,   7,   6,   5,   8,   9,   6,  10,   7,   8,   7,   8,   9.

\subsection{Firoozbakht's and related conjectures}
Farideh Firoozbakht (1962--2019) conjectured in 1982 that, if $p_n$ is the $n$th prime starting from $p_1 = 2$, 
then $(p_n)^{1/n}$ decreases strictly as $n\in\N$ increases
\cite{Rivera}.
Firoozbakht's conjecture has been verified numerically for all primes less than
$2^{64} \approx 1.844\times 10^{19}$ \cite{Visser2019}.

Ribenboim \cite[p.\ 185]{Ribenboim2004} dated this conjecture, communicated to him by the author,
 ``from about 1992,'' but 
Kourbatov \cite{Kourbatov2015} gives what appears to be the correct date, 1982 \cite{Rivera}.

Campbell and I \cite{CampbellCohen2025} first stated Conjectures \ref{conj:Fir1}--\ref{conj:Fir3}
of Firoozbakht type concerning cyclic numbers
and tested them using only the cyclics not exceeding $4\times 10^6$.
Here I test Conjectures \ref{conj:Fir1}--\ref{conj:Fir3} plus two new 
Conjectures \ref{conj:Fir4} and \ref{conj:Fir1SG} using the cyclics less than $10^8$.

\begin{conjecture}[Firoozbakht analog for cyclics 1]\label{conj:Fir1}
For every positive integer $n$ excluding $n = 1, 2, 3$ and $5$,
$$c_n^{1/n} > c_{n+1}^{1/(n+1)}.$$
The four exceptions are $1 < 2^{1/2}$,
$2^{1/2} < 3^{1/3}$,
$3^{1/3} < 5^{1/4}$, and
$7^{1/5} < 11^{1/6}$.
\end{conjecture}

\begin{conjecture}[Firoozbakht analog for cyclics 2]\label{conj:Fir2}
For every positive integer $n>1$,
$$c_n^{1/(n-1)} > c_{n+1}^{1/n}.$$
If $1^{1/0} = 1^\infty := 1$, then the only exception is $c_1 = 1 < c_2^1 = 2$.
\end{conjecture}

\begin{conjecture}[Firoozbakht analog for cyclics 3]\label{conj:Fir3}
For every $k{\in\mathbb{N}}$, there exists a least $m{\in\mathbb{N}}$, call it $N(k)$,
such that, for all $n>N(k)$,
$$c_n^{1/(n+k)} > c_{n+1}^{1/(n+k+1)}.$$
In particular, if $k = 1$ or $k = 2$, then $N(k) = 5$;
and if $k = 3$ or $k = 4$, then $N(k) = 11$.
\end{conjecture}

For the primes among the cyclic numbers, Conjecture \ref{conj:Fir1} is 
 stronger (gives tighter inequalities) than Firoozbakht's conjecture
\cite{CampbellCohen2025}. 
For example, Conjecture \ref{conj:Fir1} gives $c_6^{1/6} = 11^{1/6} \approx 
 1.4913 > c_7^{1/7} = 13^{1/7} \approx 1.4426$ whereas Firoozbakht's conjecture gives 
 $p_5^{1/5} = 11^{1/5} \approx 1.6154 > p_6^{1/6} = 13^{1/6} \approx 1.5334$. 

\begin{conjecture}[Firoozbakht analog for cyclics 4]\label{conj:Fir4}
For every $k{\in\{0\}\cup\mathbb{N}}$,
define 
$$\bar{c}(k) := \max_{n\in\N}c_n^{1/(n+k)}.$$
Then 
$$\bar{c}(0) \approx 1.4953 > \bar{c}(1) \approx 1.4085,    > \bar{c}(2)  \approx 1.3495,  >\bar{c}(3)  \approx 1.3053,  >\bar{c}(4)  \approx 1.2710>\ldots \ .$$
\end{conjecture}

\begin{conjecture}[Firoozbakht analog for SG cyclics]\label{conj:Fir1SG}
For every $n\in\N$ excluding $n = 1, 2, 3$ and $5$,
$$\sigma_n^{1/n} > \sigma_{n+1}^{1/(n+1)}.$$
The four exceptions are $1 < 2^{1/2}$,
$2^{1/2} < 3^{1/3}$,
$3^{1/3} < 5^{1/4}$, and
$7^{1/5} < 11^{1/6}$, exactly as in Conjecture \ref{conj:Fir1}.
\end{conjecture}

I confirmed Conjecture \ref{conj:Fir1SG} for the first 6,882,632 SG cyclics.

\subsection{Andrica and related conjectures}

Dorin Andrica \cite{Andrica1986} conjectured that, for all $n \in \N$,
$\Delta\sqrt{p_n} :=\sqrt{p_{n+1}}-\sqrt{p_n}< 1$.  
Visser \cite{Visser2019} verified Visser's stronger version (quoted below)
of Andrica's conjecture
for all primes less than $2^{64} \approx 1.84\times10^{19}$. 
Several claims to prove Andrica's conjecture have been published or posted
but I am not aware that any has been independently confirmed. 

\begin{conjecture}[Andrica analog for cyclics]\label{conj:Andrica}
For all $n \in \N$, $\Delta\sqrt{c_n}:=(c_{n+1})^{1/2} - (c_n)^{1/2}< 1$.
\end{conjecture}

Ribenboim \cite[p.\ 191]{Ribenboim2004} observed that the conjecture
\begin{equation}\label{eq:Ribenboim}
\lim_{n\to\infty}(\sqrt{p_{n+1}}-\sqrt{p_n}) =0
\end{equation}
would imply Andrica's conjecture. 
An editor \cite[p.\ 61]{Golomb1976} remarked that
``it is a difficult and as yet unsolved problem whether''
\eqref{eq:Ribenboim} is true.
No originator of conjecture \eqref{eq:Ribenboim} is given in \cite{Golomb1976,Ribenboim2004}.
Wolf \cite{Wolf2010} gives an impressive heuristic argument
for the truth of \eqref{eq:Ribenboim}.

\begin{figure}[htb]
\centering
\includegraphics[width=1\textwidth]{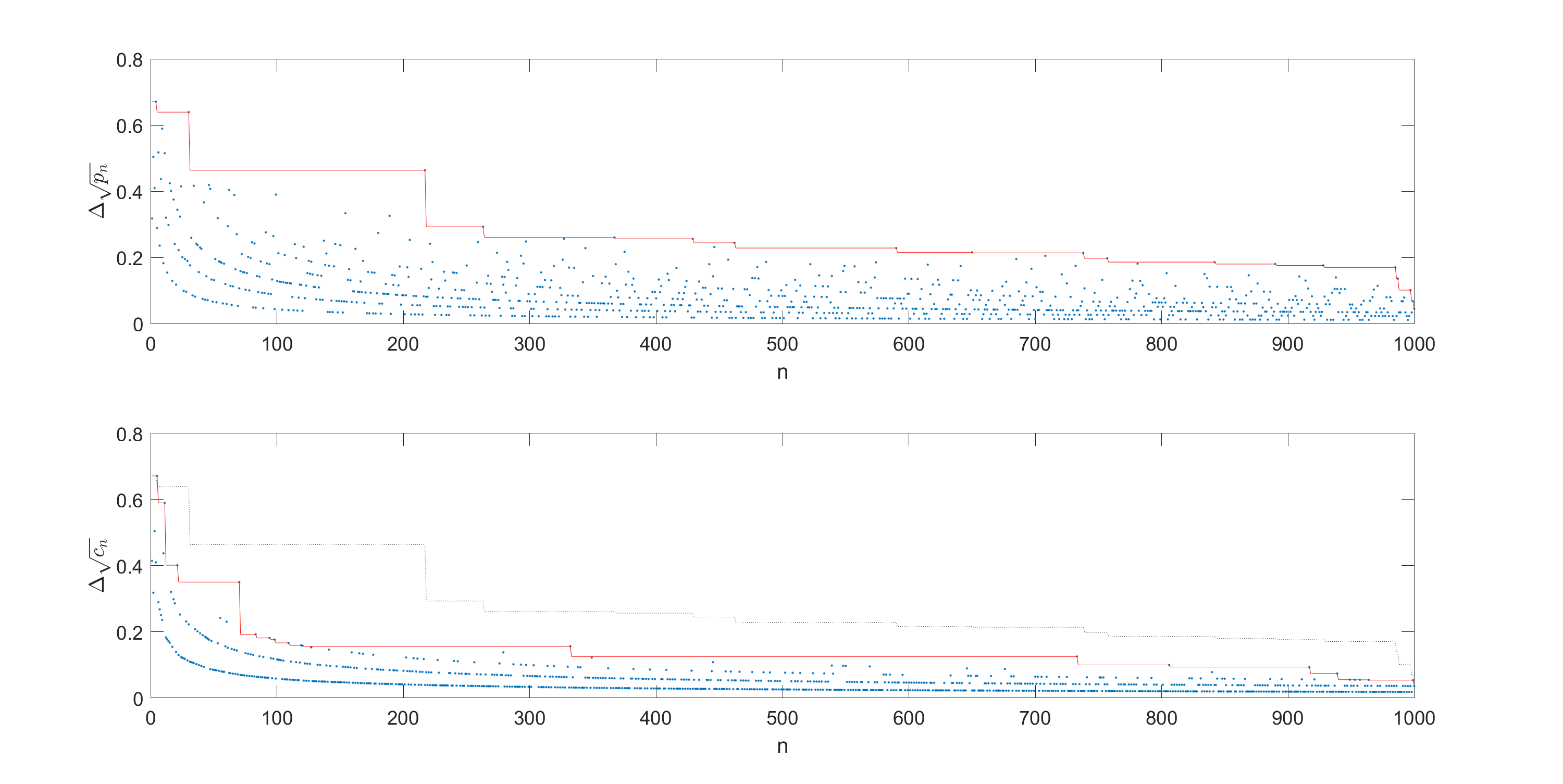}
\caption{(above) The values (blue dots) of $\Delta\sqrt{p_n} :=\sqrt{p_{n+1}}-\sqrt{p_n}$ for $n=1,\ldots,1000$ 
and the reverse cumulative maximum (red line), that is, 
the cumulative maximum starting from $\Delta\sqrt{p_{1000}}$ and working
back to $\Delta\sqrt{p_1}$. 
The decrease of the reverse cumulative maximum as $n$ increases displays the approach of $\Delta\sqrt{p_n}$
in the direction of 0 as $n$ increases.
(below) The values (blue dots) of $\Delta\sqrt{c_n} :=\sqrt{c_{n+1}}-\sqrt{c_n}$ for $n=1,\ldots,1000$ 
and the reverse cumulative maximum of $\Delta\sqrt{c_n}$ for cyclics (red line).
The grey line reproduces the reverse cumulative maximum for primes.
The grey line for primes is never less than the red line for cyclics in these examples.}
\label{fig:DeltaSqrtPrimesCyclicsCumulativeMaxReverse}
\end{figure}

Figure \ref{fig:DeltaSqrtPrimesCyclicsCumulativeMaxReverse}
provides numerical support for conjecture \eqref{eq:Ribenboim}
and for Conjecture \ref{conj:Ribenboim4cyclics}
below with $t=1/2$, which is the analogous conjecture for cyclics.
I also verified Andrica's conjecture for primes less than $10^9$ and its analog
Conjecture \ref{conj:Ribenboim4cyclics} (with $t=1/2$) for cyclics less than $10^8$.

Conjecture \ref{conj:generalizedRibenboim4primes} deals with $p_{n+1}^t-p_n^t$ for real $t\in(0,1/2]$.
Because $d^2\bigl((p+g)^t-p^t\bigr)/dg\,dt=(g + p)^t/(g + p) + t(g + p)^{t - 1}\log(g + p) > 0$ for $g>0$, $t>0$, and $p>1$,
the difference $p_{n+1}^t-p_n^t$ is an increasing function of $t$ and of $p_{n+1}-p_n$.

\begin{conjecture}[generalized analog for primes]\label{conj:generalizedRibenboim4primes}
For real $t\in(0,1/2]$,
\begin{equation}\label{eq:genRibenboimPrimes}
\lim_{n\to\infty}(p_{n+1}^t-p_n^t) =0.
\end{equation}
\end{conjecture}
Conjecture \ref{conj:generalizedRibenboim4primes} obviously implies \eqref{eq:Ribenboim}.

\begin{conjecture}[generalized analog for cyclics]\label{conj:Ribenboim4cyclics}
For real $t\in(0,1/2]$,
\begin{equation}\label{eq:genRibenboimCyclics}
\lim_{n\to\infty}(c_{n+1}^t-c_n^t) =0.
\end{equation}
\end{conjecture}

Matt Visser \cite[p.\ 182]{Visser2019strong} conjectured that, except for
$n \in \{2, 4, 6,9, 11, 30\}$ corresponding to $p_n \in \{3, 7, 13, 23, 31, 113\}$,
$$ \Delta\sqrt{p_n} :=\sqrt{p_{n+1}}-\sqrt{p_n} < \frac{1}{2}. $$
A generalization of the obvious analog for cyclics is:

\begin{conjecture}[Visser analog for cyclics]\label{conj:Visser4cyclics}
For a fraction $\epsilon\in (0, 1/2)$, there exists $N(\epsilon)\in\N$
such that, for all $n>N(\epsilon)$,
\begin{equation}\label{eq:VisserCyclics1}
\Delta\sqrt{c_n} :=\sqrt{c_{n+1}}-\sqrt{c_n} < \epsilon.
\end{equation}

In particular, except for
$n \in \{3, 5, 11\}$ corresponding to $c_n \in \{3, 7, 23\}$,
\begin{equation}\label{eq:VisserCyclics2}
\Delta\sqrt{c_n} :=\sqrt{c_{n+1}}-\sqrt{c_n} < \frac{1}{2}
\end{equation}
and except for
$n \in \{ 1,    3,    4,    5,   10,   11,   21,   70\}$ corresponding to 
$c_n \in \{1,     3,     5,     7,    19,    23,    53,   199\}$,
\begin{equation}\label{eq:VisserCyclics3}
\Delta\sqrt{c_n} :=\sqrt{c_{n+1}}-\sqrt{c_n} < \frac{1}{3}.
\end{equation}
\end{conjecture}

Alexandre Kosyak, Pieter Moree, Efthymios Sofos and Bin Zhang
(hereafter KMSZ) conjectured \cite[p.\ 216, Eq.\ (2)]{KMSZ2021}
(see also \cite{MoreeSofos2024}) that 
if $n\geq 31$ (so $p_n\geq 127$), then
$$ p_{n+1}-p_n < \sqrt{p_n}+1.$$
A generalization of this conjecture consistent with the primes less than $10^9$ is:

\begin{conjecture}[generalized KMSZ for primes]\label{conj:KMSZ4primes}
For finite positive or negative integer $k$, there exists $N(k)\in\N$ such that, for all $n>N(k)$ (strict inequality),
\begin{equation}\label{eq:genKMSZ}
p_{n+1}-p_n < \sqrt{p_n}+k.
\end{equation}
In particular, for $k\in[-20,-17]$, $N(k)=263$; for $k\in[-16,-3]$, $N(k)=217$;
for $k=-2$, $N(k)=34$;
 for $k\in[-1,+3]$, $N(k)=30$ ($k=1$ is the KMSZ conjecture);
and  for $k\geq 4$, $N(k)=0$, i.e., \eqref{eq:genKMSZ} holds for all $n\in\N$.
\end{conjecture}

An analogous conjecture is consistent with the cyclics less than $10^8$.
\begin{conjecture}[generalized KMSZ for cyclics]\label{conj:KMSZ4cyclics}
For finite positive or negative integer $k$, there exists $N(k)\in\N$ such that, for all $n>N(k)$ (strict inequality),
\begin{equation}\label{eq:genKMSZ4cyclics}
c_{n+1}-c_n < \sqrt{c_n}+k.
\end{equation}
In particular, for $k = -20, -19, -18,\ldots, -1, 0, +1$,
the corresponding 22 values of $N(k)$ are 
$N(k) = 216,   208,   176,   176,   159,   141,   127,   120, 109,   98,   83,   70,   70,   70,   70,   70,   23,   21,   21,   11,   11,   11$,
and  for $k\geq 2$, $N(k)=0$, i.e., \eqref{eq:genKMSZ4cyclics} holds for all $n\in\N$.
\end{conjecture}

As mentioned in the Introduction, Carneiro et al. \cite{Carneiro2019} proved that,
under the Riemann hypothesis, for every $p_n>3$,
$p_{n+1} - p_{n} < \frac{22}{25}\sqrt{p_n}\log p_n$.

\begin{conjecture}[Carneiro analog for cyclics]\label{conj:Carneirocyclics}
For every $c_n>3$,
\begin{equation}
c_{n+1} - c_{n} < \frac{22}{25}\sqrt{c_n}\log c_n.
\end{equation}
\end{conjecture}

The conjectured upper bounds on prime gaps in
Conjecture \ref{conj:KMSZ4primes}
and on cyclic gaps in
Conjecture \ref{conj:KMSZ4cyclics} and Conjecture \ref{conj:Carneirocyclics}
are very likely far from the best possible upper bounds.
For example,
for the primes less than $10^9$, $\max_n (p_{n+1}-p_n) = 282$ 
(the gap following prime 436273009)
while, for the next to last of these primes, 
$\sqrt{p_n}+4\approx 31626.7755$.
Similarly,
for the cyclics less than $10^8$, $\max_n (c_{n+1}-c_n) = 24$ 
while, for the next to last of these cyclics, 
$\sqrt{c_n}+2\approx 10001.9998$.
It seems worth exploring further generalizations of \eqref{eq:genKMSZ} and \eqref{eq:genKMSZ4cyclics} 
obtained by replacing the square root on the right sides by an exponent $\epsilon\in (0, 1/2)$.

\begin{conjecture}[Carneiro analog for SG cyclics]\label{conj:CarneiroSGcyclics}
For every SG cyclic $\sigma_n>3$,
\begin{equation}
\sigma_{n+1} - \sigma_{n} < \frac{22}{25}\sqrt{\sigma_n}\log \sigma_n.
\end{equation}
\end{conjecture}

I confirmed Conjecture \ref{conj:CarneiroSGcyclics} for the first 6,882,632 SG cyclics.

\subsection{Rosser, Dusart and related conjectures}

The prime number theorem \cite{Hadamard1896, delaValleePoussin189900} is equivalent to $p_n\sim n\log n$.
Rosser \cite{Rosser1939} proved that, for all $n\in\N$,
$p_n> n \log n$.
Dusart \cite{Dusart1999} proved that, for all $n >1$,
$p_n > n (\log n + \log \log n - 1).$
By analogy with Rosser's and Dusart's inequalities for primes, I conjecture, using \eqref{eq:cn} for cyclics:

\begin{conjecture}[Rosser analog for cyclics]\label{conj:Rosser}
For all $n>1$,
\begin{equation}
c_n >e^{\gamma} n\log \log \log n;\label{eq:cyclicRosser}\\
\end{equation}
\end{conjecture}

\begin{conjecture}[Dusart analog for cyclics]\label{conj:Dusart}
For all $n>1$,
\begin{equation}
c_n >e^{\gamma} n(\log \log \log n + \log \log \log \log n).\label{eq:cyclicDusart}
\end{equation}
\end{conjecture}
 Conjectures \ref{conj:Andrica},
\ref{conj:Ribenboim4cyclics},
\ref{conj:Visser4cyclics},
\ref{conj:KMSZ4cyclics},
\ref{conj:Carneirocyclics},
\ref{conj:Rosser}, and \ref{conj:Dusart}
are consistent with all numerically evaluated $c_n \in (1,10^8)$.

After seeing Conjectures \ref{conj:Rosser} and \ref{conj:Dusart} in a draft of this paper,
Carl Pomerance (personal communication, 2025-06-03 10:13) computed that
$$c_n = e^{\gamma} n (\log\log\log n  + \gamma + o(1))$$
(which refines \eqref{eq:cn}). 
So Conjecture \ref{conj:Rosser} holds for $n$ sufficiently large and 
Conjecture \ref{conj:Dusart} fails for $n$ sufficiently large.  
The numerical results in this paper fail to distinguish between 
Conjecture \ref{conj:Rosser} and Conjecture \ref{conj:Dusart} because
$\log \log \log \log 10^8 \approx 0.0670 < \gamma \approx 0.5772$.

\subsection{Additive and multiplicative inequalities}

Extending the Bertrand-Chebyshev theorem that $2p_n>p_{n+1}$ for all $n\in\N$,
Ishikawa \cite[Theorem 1]{Ishikawa1934}, Gallot et al. \cite[Lemma 10]{Gallot2011}, and I \cite[Theorem 2]{Cohen2023Bertrand} 
proved independently, and with different proofs,
that if $n>1$, then $p_n + p_{n+1} > p_{n+2}$.
When $n=1$, equality holds: $p_1 + p_2 = p_3 = 5$.

\begin{conjecture}[Ishikawa analog for cyclics]\label{conj:Ishik1}
For all $n>2$,
\begin{equation}
c_n + c_{n+1} > c_{n+2}.
\end{equation}
\end{conjecture}
 When $n=1$ or $n=2$, equality holds:  $c_1 + c_2 = c_3 = 3$ and $c_2+c_3=c_4 = 5$.

\begin{conjecture}[Ishikawa analog for SG cyclics]\label{conj:Ishik1SG}
For all $n>2$,
\begin{equation}
\sigma_n + \sigma_{n+1} > \sigma_{n+2}.
\end{equation}
\end{conjecture}
 
I verified this conjecture for the first 6,882,632 SG cyclics. 
When $n=1$ or $n=2$, equality holds:  $\sigma_1 + \sigma_2 = \sigma_3 = 3$ and $\sigma_2+\sigma_3=\sigma_4 = 5$.

More generally \cite[Theorem 1]{Cohen2023Bertrand},
if $b_1, \ldots, b_g$ are $g>1$ nonnegative integers (not necessarily distinct),
and
$d_1, \ldots, d_h$ are $h$ positive integers (not necessarily distinct),
with $1 \leq h < g$,
then there exists a positive integer $N$ 
such that, for all $n\ge N$,
$$p_{n-b_1}+p_{n-b_2}+\cdots +p_{n-b_g}>p_{n+d_1}+\cdots +p_{n+d_h}.$$
A concrete example \cite{Cohen2023Bertrand} is
$ p_{n}+p_{n+1}+p_{n+2} >  p_{n+3}+p_{n+4},$
proved for all $n\geq N=8$.
Analogously, for cyclic numbers, I conjecture:

\begin{conjecture}[sum-3-versus-sum-2 analog for cyclics]\label{conj:Cohen}
For all $n>9$,
\begin{equation}
LHS(n) := c_{n}+c_{n+1}+c_{n+2} > RHS(n) := c_{n+3}+c_{n+4}.
\end{equation}
\end{conjecture}
For $n=1, 2, 3, 4, 5, 8, 9$, I find numerically that $LHS(n)<RHS(n)$.
For $n=6, 7$ and all computed values of $n>9$, I find numerically that $LHS(n)>RHS(n)$.

According to Ribenboim \cite[p.\ 185]{Ribenboim2004}, Dusart \cite{Dusart1998thesis}
proved Mandl's conjecture that $(p_1+p_2+\cdots+p_n)/n \leq p_n/2$ for all $n>8$. 
The inequality also holds for $n=7$.
In all the examples I know, the inequality is strict wherever the weak inequality holds.

\begin{conjecture}[Dusart-Mandl analog for cyclics]\label{conj:DusartMandl}
For all $n>5$,
$$\frac{c_1+c_2+\cdots+c_n}{n} < \frac{c_n}{2}.$$
\end{conjecture}
 The opposite strict inequality holds for
$n=1, 2, 3, 4, 5$.

\begin{conjecture}[Dusart-Mandl analog for SG cyclics]\label{conj:DusartMandlSG}
For all $n>5$,
$$\frac{\sigma_1+\sigma_2+\cdots+\sigma_n}{n} < \frac{\sigma_n}{2}.$$
\end{conjecture}
 I verified this conjecture for the first 6,882,632 SG cyclics. 
The opposite strict inequality holds for
$n=1, 2, 3, 4, 5$.

Lauren\c{t}iu Panaitopol (1940--2008) \cite{Panaitopol1998} and I
\cite{Cohen2022} proved independently a multiplicative inequality for primes:
if $m\in\N,\ n\in\N$, and $m\leq n$, then $p_{m\cdot n} < p_mp_n$
unless $(m,n) = (3,4)$ or $(m,n) = (4,4)$, in which cases
the reverse strict inequality holds.

\begin{conjecture}[Panaitopol analog for cyclics]\label{conj:Panaitopol}
If $m\in\N,\ n\in\N$, and $3\leq m\leq n$, then $c_{m\cdot n} < c_mc_n$
unless $(m,n) = (3,3)$ or  $(m,n) = (5,h)$ for $h = 5, 6, 7, 8, 9, 10$, in which cases
the reverse strict inequality holds.
\end{conjecture}

Anton Vrba \cite{Vrba} conjectured (in 2010, according to Kourbatov \cite{Kourbatov2016}) that
\begin{equation}\label{pnovergeomean}
\lim_{n\to\infty} \frac{p_n}{\left(\prod_{j=1}^n p_j\right)^{1/n}} = e.
\end{equation}
S\'andor and Verroken \cite{SandorVerroken2011} in 2011 and independently Kourbatov \cite{Kourbatov2016} in 2016 proved \eqref{pnovergeomean}.
Hassani \cite{Hassani2013, Hassani2014} and Kourbatov \cite{Kourbatov2016} bounded the approach of
$p_n/(\prod_{j=1}^n p_j)^{1/n}$ to $e$.
Kourbatov \cite{Kourbatov2016} gave a short proof and calculated higher-order terms in a series expansion.
\begin{conjecture}[Vrba analog for cyclics]\label{conj:Vrba}
$$ \lim_{n\to\infty} \frac{c_n}{\left(\prod_{j=1}^n c_j\right)^{1/n}} = e.$$
\end{conjecture}
 For $n=28488167$, I calculate $c_n = 99999997$ and ${c_n}/{(\prod_{j=1}^n c_j)^{1/n}} \approx2.7362$,
not a bad approximation to $e \approx 2.7183$. 
To circumvent numerical overflow of $\prod_{j=1}^n c_j$ in this calculation, 
I computed
${c_n}/{(\prod_{j=1}^n c_j)^{1/n}}$ by means of the equivalent $\exp\{\log c_n - (\sum_{j=1}^n \log  c_j)/n\}$.

Hassani \cite[Theorem 1.1]{Hassani2013} also proved that 
$$\lim_{n\to\infty} \frac{(p_1 + \cdots + p_n)/n}{(p_1 \times \cdots \times p_n)^{1/n}} = \frac{e}{2}.$$

\begin{conjecture}[Hassani analog for cyclics]\label{conj:Hassani}
$$\lim_{n\to\infty} \frac{(c_1 + \cdots + c_n)/n}{(c_1 \times \cdots \times c_n)^{1/n}} = \frac{e}{2}.$$
\end{conjecture}
For $n=28488167$, I calculate
$$\frac{(c_1 + \cdots + c_n)/n}{(c_1 \times \cdots \times c_n)^{1/n}} \approx 1.3638,\qquad \frac{e}{2}\approx 1.3591.$$

If Conjectures \ref{conj:Vrba} and \ref{conj:Hassani} are true, then
dividing the former equality by the latter equality yields
an equality Campbell and I \cite[Theorem 1]{CampbellCohen2025} proved:
$$ \lim_{n\to\infty} \frac{c_n}{(c_1 + \cdots + c_n)/n} = {2}. $$
Thus the difference between $(c_1+c_2+\cdots+c_n)/n$ and $c_n/2$,
which appear in an inequality
in Conjecture \ref{conj:DusartMandl},
is proved to vanish asymptotically as $n\to\infty$.

Analogs of Conjectures \ref{conj:Vrba} and \ref{conj:Hassani} 
for SG cyclics, replacing $c_n$ by $\sigma_n$, are obvious
and are numerically plausible.

\subsection{Sequence of absolute difference sequences: Proth, Gilbreath}

Define $\{a(n) \mid n\in \{0\} \cup \N\}$ to be a \emph{sequence of absolute difference 
sequences} (hereafter, SADS) if, for each 
$n\in \{0\} \cup \N$, $a(n) = (a_1(n), a_2(n),a_3(n),\ldots)$ is an infinite sequence of real numbers
such that, for all $n\in\N$ and all $m\in\N$, $a_m(n) = |a_m(n-1)-a_{m+1}(n-1)|$. 
In more detail, starting from an arbitrary initial real sequence 
$a(0):=(a_1(0), a_2(0),a_3(0),\ldots,a_m(0),a_{m+1}(0),\ldots)$,
the next sequence is
$a(1):=(|a_1(0)- a_2(0)|$, $|a_2(0)-a_3(0)|,\ldots,|a_m(0)-a_{m+1}(0)|,\ldots)$,
and the elements of each successive sequence $a(n)$ give the absolute values of the first differences of
elements of the preceding sequence.

Fran\c{c}ois Proth (1852--1879) \cite{Proth1878} in 1878 computed numerically the behavior of a SADS starting
from the first seven primes (he included 1 as the first prime, which I ignore) and observed that
every successor sequence begins with 1. For example, $a_1(1) = |3-2| = 1$ 
and $a_1(2) = | |3-2| - |5-3| | = |1 - 2| = 1$, and so on.
He suggested that $a_1(n) = 1$ for all $n\in\N$.
The journal editor, Eug\`ene Catalan, in a gentle concluding footnote, asked (my translation): 
``Are not the \emph{theorems} of M. Proth that one has just read rather postulates?'' 
[$\ll$Est-ce que les \emph{th\'eor\`emes} de M. Proth, 
qu'on vient de lire, ne sont pas, plut\^{o}t, des postulata?$\gg$]
I take Proth's suggestion as a conjecture, generally known as N. L. Gilbreath's conjecture
\cite[pp.\ 191-192]{Ribenboim2004}.
The Proth-Gilbreath conjecture has been verified for  the primes less than $10^{13}$.

\begin{conjecture}[Proth-Gilbreath analog for cyclics]\label{conj:Proth}
A SADS starting
from the cyclic numbers $\C$ after omitting $c_1=1$ has 1
as the first element of every successor sequence.
\end{conjecture}
I verified this conjecture for 1 million successor sequences of the cyclic numbers following but not including $c_1=1$.

\begin{conjecture}[Proth-Gilbreath analog for SG cyclics]\label{conj:ProthSG}
A SADS starting
from the SG cyclics after omitting $\sigma_1=1$ has 1
as the first element of every successor sequence.
\end{conjecture}
I verified this conjecture for 1 million successor sequences of the SG cyclics following but not including $\sigma_1=1$.

\section{Second Hardy and Littlewood conjecture: cyclic analog is false}\label{sec:HLconjecture2}

The second conjecture of Hardy and Littlewood 
\cite{HardyLittlewood1923}
states that $\pi(m+n)\leq \pi(m) + \pi(n)$ 
for all integers $2\leq m \leq n$.

\begin{conjecture}[Hardy and Littlewood analog for SG primes]\label{conj:HardyLittlewood2SGprimes}
The counting function of SG primes obeys
$\pi_{SG}(m+n)\leq \pi_{SG}(m) + \pi_{SG}(n)$ for all $2\leq m \leq n$.
\end{conjecture}

I verified Conjecture \ref{conj:HardyLittlewood2SGprimes} for all
$2\leq m \leq 910664$, $m \leq n \leq 999997$.

Using the
counting function $C(\cdot)$ \eqref{eq:cycliccountingfn} of cyclic numbers,
an analog for cyclics of the second conjecture of Hardy and Littlewood, 
starting from $1\leq m \leq n$, is:
\begin{conjecture}[Hardy and Littlewood analog for cyclics]\label{conj:HardyLittlewood}
For all integers $1\leq m \leq n$, $C(m+n)\leq C(m) + C(n)$.
\end{conjecture}
This conjecture is false. Let $m  = 209 = c_{71}$, so $C(209) = 71$.
Let $n = 389 =  c_{128}$, so $C(389) = 128$. Then
$C(m+n) = C(598) = 200$ because $c_{200}= 595<598<c_{201}=599$.
Thus  $C(m+n) = C(598) = 200> C(m) + C(n) = 71 + 128 = 199$.
Counterexamples like this one are abundant.

To the extent that analogies between primes and cyclics are valid, 
this counterexample gives a further small hint to support the belief of 
Hensley and Richards \cite{HensleyRichards1974} 
that the second conjecture of Hardy and Littlewood \cite{HardyLittlewood1923} about primes is  false.

\begin{conjecture}[Hardy and Littlewood analog for Sophie Germain cyclics]\label{conj:HardyLittlewoodSG}
For all integers $1\leq m \leq n$, $C_{\sigma}(m+n)\leq C_{\sigma}(m) + C_{\sigma}(n)$.
\end{conjecture}

The counterexample to Conjecture \ref{conj:HardyLittlewood}
for cyclics, $m  = 209 = \sigma_{46}$, 
$n = 389 =  \sigma_{83}$ is not a counterexample
to Conjecture \ref{conj:HardyLittlewoodSG} for SG cyclics because 
$\sigma_{120}= 593<598<\sigma_{121}=599$ and therefore
$C_{\sigma}(209+389) = C_{\sigma}(598) = 120 <$
$C_{\sigma}(209) + C_{\sigma}(389) =  46+ 83= 129$. 

For $m=1,\ldots,10^6$ and $n=m,\ldots, 10^6$, I found no counterexamples to Conjecture \ref{conj:HardyLittlewoodSG} for SG cyclics. It remains open.

\section{Acknowledgments}
My enthusiastic thanks to
John M. Campbell for introducing me to cyclic numbers and for our collaboration 
\cite{CampbellCohen2025};
Barry Mazur (personal communication, 2024-11-17 18:08) for proposing the essential idea of sequence $B(\cdot)$, closely related to the sequence $L(\cdot)$;
Pierre Deligne (personal communication, 2025-03-06 19:12) for improving a conjecture in
a previous draft with his Conjecture \ref{conj:Np};
Carl Pomerance (personal communications, 2025-03-10 14:27 and 17:49) \cite{Pomerance2025}, after seeing a previous draft,
for proving, improving, partially proving, or disproving Conjectures \ref{conj:Goldbach}, \ref{conj:twincyclics}, \ref{conj:cyclictriplets}, 
\ref{conj:GermainCyclics}, \ref{conj:SGcyclicsMod3}, \ref{conj:Rosser}, and \ref{conj:Dusart}; Alexei Kourbatov for
a citation \cite{Rivera} to the correct date of Firoozbakht's conjecture and suggesting
evaluation of Lagneau's condition by modular exponentiation \cite{HAC1996, WikiModularExponentiation};
Alexei Kourbatov, Pieter Moree, Elchin Hasanalizade, and Marek Wolf
for making me aware of relevant works \cite{Kourbatov2015,Gallot2011,KMSZ2021,Panaitopol1998, MoreeSofos2024,
Wolf1998,Wolf2010};
Richard P. Stanley for permission to quote his proof that the conditions of Szele \cite{Szele1947}
and Lagneau \cite[\seqnum{A003277}, November 18, 2012]{oeis} are equivalent;
Chris Higgins, Eric Weisstein, and the Wolfram content team for providing
the exact citation \cite{brocard1904} and access to Brocard's conjecture;
Mehdi Hassani for an informative exchange;
Jeffrey O. Shallit (personal communication, 2025-01-01 14:41)
for finding and sending me Golubew's article \cite{Golubew1957};
Adam C. Green (personal communication, 2025-05-29 03:42), Archivist, 
Trinity College Library, Cambridge, for information about Alfred Edward Western;
and Roseanne Benjamin for help during this work. 
I thank the reviewer of a previous draft for several helpful comments,
including a question that led to Theorem \ref{thm:sumofcyclics}.

\bigskip
\hrule
\bigskip

\noindent 2020 {\it Mathematics Subject Classification}:
Primary 11N25; Secondary 11N05, 11N37, 11N56, 11N60, 11N69.

\noindent \emph{Keywords: } cyclic number, prime number, Landau's problems.

\bigskip
\hrule
\bigskip

\noindent (Concerned with sequences
\seqnum{A000040}, \seqnum{A000720}, \seqnum{A002496}, \seqnum{A003277}, \seqnum{A014085}, \seqnum{A050216}, \seqnum{A050384}, \seqnum{A061091}, and \seqnum{A349997}.)

\bigskip
\hrule
\bigskip

\vspace*{+.1in}
\noindent
Received January 5, 2025;
revised versions received  June 12, 2025, .
Published in {\it Journal of Integer Sequences},

\bigskip
\hrule
\bigskip

\noindent
Return to
{\em Journal of Integer Sequences} home page \url{https://cs.uwaterloo.ca/journals/JIS/}
\vskip .1in

\end{document}